\documentclass[reqno,12pt,a4paper]{amsart}

\usepackage{mathrsfs}
\usepackage{amssymb}
\usepackage{amsfonts}
\usepackage{latexsym}
\usepackage{amsthm}
\usepackage{graphicx}
\usepackage{xcolor}
\def\lb{\label}

\newcommand{\er}[1]{\textrm{(\ref{#1})}}

\begin{document}

\renewcommand{\theequation}{\arabic{section}.\arabic{equation}}
\theoremstyle{plain}
\newtheorem{theorem}{\bf Theorem}[section]
\newtheorem{lemma}[theorem]{\bf Lemma}
\newtheorem{corollary}[theorem]{\bf Corollary}
\newtheorem{proposition}[theorem]{\bf Proposition}
\newtheorem{definition}[theorem]{\bf Definition}

\theoremstyle{remark}
\newtheorem{remark}[theorem]{\bf Remark}
\newtheorem{example}[theorem]{\bf Example}

\def\a{\alpha}  \def\cA{{\mathcal A}}     \def\bA{{\bf A}}  \def\mA{{\mathscr A}}
\def\b{\beta}   \def\cB{{\mathcal B}}     \def\bB{{\bf B}}  \def\mB{{\mathscr B}}
\def\g{\gamma}  \def\cC{{\mathcal C}}     \def\bC{{\bf C}}  \def\mC{{\mathscr C}}
\def\G{\Gamma}  \def\cD{{\mathcal D}}     \def\bD{{\bf D}}  \def\mD{{\mathscr D}}
\def\d{\delta}  \def\cE{{\mathcal E}}     \def\bE{{\bf E}}  \def\mE{{\mathscr E}}
\def\D{\Delta}  \def\cF{{\mathcal F}}     \def\bF{{\bf F}}  \def\mF{{\mathscr F}}
\def\c{\chi}    \def\cG{{\mathcal G}}     \def\bG{{\bf G}}  \def\mG{{\mathscr G}}
\def\z{\zeta}   \def\cH{{\mathcal H}}     \def\bH{{\bf H}}  \def\mH{{\mathscr H}}
\def\e{\eta}    \def\cI{{\mathcal I}}     \def\bI{{\bf I}}  \def\mI{{\mathscr I}}
\def\p{\psi}    \def\cJ{{\mathcal J}}     \def\bJ{{\bf J}}  \def\mJ{{\mathscr J}}
\def\vT{\Theta} \def\cK{{\mathcal K}}     \def\bK{{\bf K}}  \def\mK{{\mathscr K}}
\def\k{\kappa}  \def\cL{{\mathcal L}}     \def\bL{{\bf L}}  \def\mL{{\mathscr L}}
\def\l{\lambda} \def\cM{{\mathcal M}}     \def\bM{{\bf M}}  \def\mM{{\mathscr M}}
\def\L{\Lambda} \def\cN{{\mathcal N}}     \def\bN{{\bf N}}  \def\mN{{\mathscr N}}
\def\m{\mu}     \def\cO{{\mathcal O}}     \def\bO{{\bf O}}  \def\mO{{\mathscr O}}
\def\n{\nu}     \def\cP{{\mathcal P}}     \def\bP{{\bf P}}  \def\mP{{\mathscr P}}
\def\r{\rho}    \def\cQ{{\mathcal Q}}     \def\bQ{{\bf Q}}  \def\mQ{{\mathscr Q}}
\def\s{\sigma}  \def\cR{{\mathcal R}}     \def\bR{{\bf R}}  \def\mR{{\mathscr R}}
\def\S{\Sigma}  \def\cS{{\mathcal S}}     \def\bS{{\bf S}}  \def\mS{{\mathscr S}}
\def\t{\tau}    \def\cT{{\mathcal T}}     \def\bT{{\bf T}}  \def\mT{{\mathscr T}}
\def\f{\phi}    \def\cU{{\mathcal U}}     \def\bU{{\bf U}}  \def\mU{{\mathscr U}}
\def\F{\Phi}    \def\cV{{\mathcal V}}     \def\bV{{\bf V}}  \def\mV{{\mathscr V}}
\def\P{\Psi}    \def\cW{{\mathcal W}}     \def\bW{{\bf W}}  \def\mW{{\mathscr W}}
\def\o{\omega}  \def\cX{{\mathcal X}}     \def\bX{{\bf X}}  \def\mX{{\mathscr X}}
\def\x{\xi}     \def\cY{{\mathcal Y}}     \def\bY{{\bf Y}}  \def\mY{{\mathscr Y}}
\def\X{\Xi}     \def\cZ{{\mathcal Z}}     \def\bZ{{\bf Z}}  \def\mZ{{\mathscr Z}}
\def\be{{\bf e}} \def\bc{{\bf c}}
\def\bv{{\bf v}} \def\bu{{\bf u}}
\def\Om{\Omega}
\def\bp{{\bf p}}\def\bq{{\bf q}}
\def\bx{{\bf x}} \def\by{{\bf y}}
\def\mm{\mathrm m}
\def\mn{\mathrm n}
\def\rmH{\mathrm H}
\newcommand{\mc}{\mathscr {c}}
\newcommand{\gA}{\mathfrak{A}}          \newcommand{\ga}{\mathfrak{a}}
\newcommand{\gB}{\mathfrak{B}}          \newcommand{\gb}{\mathfrak{b}}
\newcommand{\gC}{\mathfrak{C}}          \newcommand{\gc}{\mathfrak{c}}
\newcommand{\gD}{\mathfrak{D}}          \newcommand{\gd}{\mathfrak{d}}
\newcommand{\gE}{\mathfrak{E}}
\newcommand{\gF}{\mathfrak{F}}           \newcommand{\gf}{\mathfrak{f}}
\newcommand{\gG}{\mathfrak{G}}           
\newcommand{\gH}{\mathfrak{H}}           \newcommand{\gh}{\mathfrak{h}}
\newcommand{\gI}{\mathfrak{I}}           \newcommand{\gi}{\mathfrak{i}}
\newcommand{\gJ}{\mathfrak{J}}           \newcommand{\gj}{\mathfrak{j}}
\newcommand{\gK}{\mathfrak{K}}            \newcommand{\gk}{\mathfrak{k}}
\newcommand{\gL}{\mathfrak{L}}            \newcommand{\gl}{\mathfrak{l}}
\newcommand{\gM}{\mathfrak{M}}            \newcommand{\gm}{\mathfrak{m}}
\newcommand{\gN}{\mathfrak{N}}            \newcommand{\gn}{\mathfrak{n}}
\newcommand{\gO}{\mathfrak{O}}
\newcommand{\gP}{\mathfrak{P}}             \newcommand{\gp}{\mathfrak{p}}
\newcommand{\gQ}{\mathfrak{Q}}             \newcommand{\gq}{\mathfrak{q}}
\newcommand{\gR}{\mathfrak{R}}             \newcommand{\gr}{\mathfrak{r}}
\newcommand{\gS}{\mathfrak{S}}              \newcommand{\gs}{\mathfrak{s}}
\newcommand{\gT}{\mathfrak{T}}             \newcommand{\gt}{\mathfrak{t}}
\newcommand{\gU}{\mathfrak{U}}             \newcommand{\gu}{\mathfrak{u}}
\newcommand{\gV}{\mathfrak{V}}             \newcommand{\gv}{\mathfrak{v}}
\newcommand{\gW}{\mathfrak{W}}             \newcommand{\gw}{\mathfrak{w}}
\newcommand{\gX}{\mathfrak{X}}               \newcommand{\gx}{\mathfrak{x}}
\newcommand{\gY}{\mathfrak{Y}}              \newcommand{\gy}{\mathfrak{y}}
\newcommand{\gZ}{\mathfrak{Z}}             \newcommand{\gz}{\mathfrak{z}}

\def\ve{\varepsilon}   \def\vt{\vartheta}    \def\vp{\varphi}    \def\vk{\varkappa} \def\vr{\varrho}

\def\A{{\mathbb A}} \def\B{{\mathbb B}} \def\C{{\mathbb C}}
\def\dD{{\mathbb D}} \def\E{{\mathbb E}} \def\dF{{\mathbb F}} \def\dG{{\mathbb G}} \def\H{{\mathbb H}}\def\I{{\mathbb I}} \def\J{{\mathbb J}} \def\K{{\mathbb K}} \def\dL{{\mathbb L}}\def\M{{\mathbb M}} \def\N{{\mathbb N}} \def\O{{\mathbb O}} \def\dP{{\mathbb P}} \def\R{{\mathbb R}}\def\S{{\mathbb S}} \def\T{{\mathbb T}} \def\U{{\mathbb U}} \def\V{{\mathbb V}}\def\W{{\mathbb W}} \def\X{{\mathbb X}} \def\Y{{\mathbb Y}} \def\Z{{\mathbb Z}}


\def\la{\leftarrow}              \def\ra{\rightarrow}            \def\Ra{\Rightarrow}
\def\ua{\uparrow}                \def\da{\downarrow}
\def\lra{\leftrightarrow}        \def\Lra{\Leftrightarrow}


\def\lt{\biggl}                  \def\rt{\biggr}
\def\ol{\overline}               \def\wt{\widetilde}
\def\ul{\underline}
\def\no{\noindent}


\let\ge\geqslant
\let\le\leqslant
\def\lan{\langle}                \def\ran{\rangle}
\def\/{\over}                    \def\iy{\infty}
\def\sm{\setminus}               \def\es{\emptyset}
\def\ss{\subset}                 \def\ts{\times}
\def\pa{\partial}                \def\os{\oplus}
\def\om{\ominus}                 \def\ev{\equiv}
\def\iint{\int\!\!\!\int}        \def\iintt{\mathop{\int\!\!\int\!\!\dots\!\!\int}\limits}
\def\el2{\ell^{\,2}}             \def\1{1\!\!1}
\def\sh{\sharp}
\def\wh{\widehat}
\def\bs{\backslash}
\def\intl{\int\limits}


\def\na{\mathop{\mathrm{\nabla}}\nolimits}
\def\sh{\mathop{\mathrm{sh}}\nolimits}
\def\ch{\mathop{\mathrm{ch}}\nolimits}
\def\where{\mathop{\mathrm{where}}\nolimits}
\def\all{\mathop{\mathrm{all}}\nolimits}
\def\as{\mathop{\mathrm{as}}\nolimits}
\def\Area{\mathop{\mathrm{Area}}\nolimits}
\def\arg{\mathop{\mathrm{arg}}\nolimits}
\def\const{\mathop{\mathrm{const}}\nolimits}
\def\det{\mathop{\mathrm{det}}\nolimits}
\def\diag{\mathop{\mathrm{diag}}\nolimits}
\def\diam{\mathop{\mathrm{diam}}\nolimits}
\def\dim{\mathop{\mathrm{dim}}\nolimits}
\def\dist{\mathop{\mathrm{dist}}\nolimits}
\def\Im{\mathop{\mathrm{Im}}\nolimits}
\def\Iso{\mathop{\mathrm{Iso}}\nolimits}
\def\Ker{\mathop{\mathrm{Ker}}\nolimits}
\def\Lip{\mathop{\mathrm{Lip}}\nolimits}
\def\rank{\mathop{\mathrm{rank}}\limits}
\def\Ran{\mathop{\mathrm{Ran}}\nolimits}
\def\Re{\mathop{\mathrm{Re}}\nolimits}
\def\Res{\mathop{\mathrm{Res}}\nolimits}
\def\res{\mathop{\mathrm{res}}\limits}
\def\sign{\mathop{\mathrm{sign}}\nolimits}
\def\supp{\mathop{\mathrm{supp}}\nolimits}
\def\Tr{\mathop{\mathrm{Tr}}\nolimits}
\def\BBox{\hspace{1mm}\vrule height6pt width5.5pt depth0pt \hspace{6pt}}


\newcommand\nh[2]{\widehat{#1}\vphantom{#1}^{(#2)}}
\def\dia{\diamond}

\def\Oplus{\bigoplus\nolimits}



\def\qqq{\qquad}
\def\qq{\quad}
\let\ge\geqslant
\let\le\leqslant
\let\geq\geqslant
\let\leq\leqslant
\newcommand{\ca}{\begin{cases}}
\newcommand{\ac}{\end{cases}}
\newcommand{\ma}{\begin{pmatrix}}
\newcommand{\am}{\end{pmatrix}}
\renewcommand{\[}{\begin{equation}}
\renewcommand{\]}{\end{equation}}
\def\eq{\begin{equation}}
\def\qe{\end{equation}}
\def\[{\begin{equation}}
\def\bu{\bullet}

\title
{Asymptotic isospectrality of Schr\"odinger operators on periodic graphs}

\date{\today}
\author[Natalia Saburova]{Natalia Saburova}
\address{Northern (Arctic) Federal University, Severnaya Dvina Emb. 17, Arkhangelsk, 163002, Russia,
 \ n.saburova@gmail.com, \ n.saburova@narfu.ru}

\subjclass{} \keywords{discrete Schr\"odinger operators, periodic graphs, adding an edge, asymptotics of spectral band edges, isospectrality}

\begin{abstract}
We consider discrete Schr\"odinger operators with periodic potentials on periodic graphs. Their spectra consist of a finite number of bands. We perturb a periodic graph by adding edges in a periodic way (without changing the vertex set) and show that if the added edges are long enough, then the perturbed graph is asymptotically isospectral to some periodic graph of a higher dimension but without long edges. We also obtain a criterion for the perturbed graph to be not only asymptotically isospectral but just isospectral to this higher dimensional periodic graph. One of the simplest examples of such asymptotically isospectral periodic graphs is the square lattice perturbed by long edges and the cubic lattice. We also get asymptotics of the endpoints of the spectral bands for the Schr\"odinger operator on the perturbed graph as the length of the added edges tends to infinity.
\end{abstract}

\maketitle

\section {\lb{Sec0}{Introduction}}
\setcounter{equation}{0}
The spectral analysis of Schr\"odinger operators on periodic structures is one of the most important problems in solid state physics. Spectral properties of these operators determine main physical properties of the material. It is well-known that the spectrum of periodic Schr\"odinger operators is the union of closed finite intervals called \emph{spectral bands} \cite{RS78}. The spectral bands may leave open intervals between them called \emph{spectral gaps} (forbidden zones). The spectral bands correspond to the set of wave energies at which the waves can propagate through the medium. A spectral gap describes a set of wave-lengths for which no transport is permitted through the media. The band edges mark the boundary between propagation and insulation. For example in semiconductors, an electron with an energy in the forbidden zone cannot propagate through the material. Similarly, in a photonic crystal, an optical counterpart of semi-conductors, a spectral gap shows that electromagnetic waves with certain frequencies cannot propagate through the photonic crystal (see, e.g., \cite{KK02}).

In this paper we consider discrete Schr\"odinger operators on periodic graphs. Physically, vertices of a graph can be thought as a discretization of space and the Laplacian on $\R^d$ is replaced by a finite difference operator. The vertices can also be seen as locations of atoms in a solid. In this atomic interpretation the edges become chemical bonding of atoms, and the model is known as the tight-binding approximation \cite{AM76}.

The spectrum of the discrete Schr\"odinger operators with periodic potentials on periodic graphs consists of a finite number of bands. The aim  of this paper is twofold. Firstly, we study behavior of the spectrum when edges are added to a periodic graph in a periodic way (without changing the vertex set). Secondly, we give a general construction of non-isomorphic isospectral periodic graphs and periodic graphs with asymptotically close spectra. By isospectrality of periodic graphs we mean that the spectra of the Schr\"odinger operators on the graphs coincide as sets.

Discrete Schr\"odinger operators for a class of periodic graphs with compact perturbations including the square, triangular, diamond, kagome lattices were discussed in \cite{AIM16}. Laplacians on periodic graphs with non-compact perturbations and the stability of their essential spectrum were considered in \cite{SaSu17}. The spectrum of Laplacians on the lattice $\Z^d$ with pendant edges was studied in \cite{Su13}. In \cite{KS17} the authors considered Laplace operators on periodic graphs perturbed by \emph{guides}, i.e., graphs that are periodic in some directions and finite in other ones. A procedure of creating gaps in the spectrum of Laplacians on an infinite graph $\cG$ by attaching a fixed finite graph $\cG_o$ to each vertex of $\cG$ was presented in \cite{AS00}. Spectral and scattering theory for discrete Schr\"odinger operators on periodic graphs perturbed by adding infinitely many edges with some suitable weights on them and/or by removing a finite number of edges was  investigated in \cite{RT22}. In that case, in general, the perturbed graph is not periodic anymore. For finite graphs, effect of deleting edges on the spectrum of discrete magnetic Laplacians was systematically studied in \cite{FLP20} (see also the references therein).

\medskip

In this paper we consider perturbations of a periodic graph by adding edges between existing vertices in a periodic way. We describe our main results:

$\bullet$ we show that if the added edges are long enough, then the perturbed graph is asymptotically isospectral to some periodic graph $\wt\cG$ of a higher dimension but without long edges (i.e., edges connecting distant vertices);

$\bullet$ we obtain asymptotics of the endpoints of the spectral bands for the Schr\"odinger operator on the perturbed graph as the length of the added edges tends to infinity;

$\bullet$ we formulate a criterion for the perturbed graph to be not only asymptotically isospectral but just isospectral to the higher dimensional periodic graph $\wt\cG$.

\medskip

We note that the spectrum of the discrete Laplacian on \emph{finite graphs} consists of a finite number of eigenvalues. Finite graphs are isospectral if they have the same eigenvalues with the same multiplicity. Construction techniques of isospectral finite graphs are presented in \cite{BPBS09,FLP23} (see also the references therein).

The spectrum of the operators on \emph{periodic graphs} has a much more complicated band-gap structure and there are some notions of isospectrality for periodic graphs: Floquet isospectrality, Fermi isospectrality, isospectrality for the periodic problem (see Section \ref{Siso} for the definitions and references). In this paper by isospectrality of periodic graphs we just mean that the spectra of the Schr\"odinger operators on the graphs coincide as sets.

\subsection{Periodic graphs and edge indices.}
Let $\cG=(\cV,\cE)$ be a connected infinite graph, possibly  having
loops and multiple edges and embedded into the space $\R^d$. Here
$\cV$  is the set of its vertices and $\cE$ is the set of its
unoriented edges. If a vertex is an endpoint of an edge, they are  \emph{incident}. Considering each edge in $\cE$ to have two
orientations, we introduce the set $\cA$ of all oriented edges. An
edge starting at a vertex $u$ and ending at a vertex $v$ from $\cV$
will be denoted as the ordered pair $(u,v)\in\cA$. Let $\ol\be=(v,u)$ be the inverse edge of $\be=(u,v)\in\cA$. We define the \emph{degree} ${\vk}_v$ of a vertex $v\in\cV$ as the number of all unoriented edges from $\cE$ incident to the vertex $v$ (with loops counted twice),  or, equivalently, as the number of all oriented edges from $\cA$ starting at $v$.

Let $\G$ be a lattice of rank $d$ in $\R^d$ with a basis $\{\ga_1,\ldots,\ga_d\}$, i.e.,
$$
\textstyle\G=\Big\{\ga\in\R^d : \ga=\sum\limits_{s=1}^dn_s\ga_s, \; (n_s)_{s=1}^d\in\Z^d\Big\},
$$
and let
\[\lb{fuce}
\textstyle\Omega=\Big\{\bx\in\R^d : \bx=\sum\limits_{s=1}^dx_s\ga_s, \; (x_s)_{s=1}^d\in
[0,1)^d\Big\}
\]
be the \emph{fundamental cell} of the lattice $\G$. We define the equivalence relation on $\R^d$:
$$
\bx\equiv \by \; (\hspace{-4mm}\mod \G) \qq\Leftrightarrow\qq \bx-\by\in\G \qqq
\forall\, \bx,\by\in\R^d.
$$

We consider \emph{locally finite $\G$-periodic graphs} $\cG$, i.e.,
graphs satisfying the following conditions:
\begin{itemize}
  \item[1)] $\cG=\cG+\ga$ for any $\ga\in\G$, i.e., $\cG$ is invariant under translation by any vector $\ga\in\G$;
  \item[2)] the quotient graph  $\cG_*=\cG/\G$ is finite.
\end{itemize}
The basis $\ga_1,\ldots,\ga_d$ of the lattice $\G$ is called the {\it periods}  of $\cG$. We also call the quotient graph $\cG_*=\cG/\G$ the \emph{fundamental graph} of the periodic graph $\cG$. The fundamental graph $\cG_*=(\cV_*,\cE_*)$ has the vertex set $\cV_*=\cV/\G$, the set $\cE_*=\cE/\G$ of unoriented edges and the set $\cA_*=\cA/\G$ of oriented edges which are finite.

\medskip

We define the important notion of an {\it edge index} which was introduced in \cite{KS14}. This notion allows one to consider, instead of a periodic graph, the finite fundamental graph with edges labeled by some integer vectors called indices.

For each $\bx\in\R^d$, we introduce the vector $\bx_\A\in\R^d$ by
\[\lb{cola}
\bx_\A=(x_1,\ldots,x_d), \qqq \textrm{where} \qq \bx=
\textstyle\sum\limits_{s=1}^dx_s\ga_s,
\]
i.e., $\bx_\A$ is the coordinate vector of $\bx$ with respect to the
basis  $\A=\{\ga_1,\ldots,\ga_d\}$ of the lattice~$\G$.

For any vertex $v\in\cV$ of a $\G$-periodic graph $\cG$, the
following  unique representation holds true:
\[\lb{Dv}
v=v_0+[v], \qqq \textrm{where}\qqq v_0\in\cV\cap\Omega,\qqq [v]\in\G,
\]
and $\Omega$ is the fundamental cell of the lattice $\G$ defined by
\er{fuce}. In other words, each vertex $v$ can be obtained from a
vertex $v_0\in \Omega$  by a shift by a vector $[v]\in\G$.

For any oriented edge $\be=(u,v)\in\cA$ of the $\G$-periodic graph
$\cG$, we define the \emph{edge index} $\t(\be)$ as the vector of the lattice $\Z^d$ given by
\[
\lb{in}
\t(\be)=[v]_\A-[u]_\A\in\Z^d,
\]
where $[v]\in\G$ is defined by \er{Dv} and the vector $[v]_\A\in\Z^d$  is given by \er{cola}. Due to periodicity of the graph $\cG$, the edge indices $\t(\be)$ satisfy
\[\lb{Gpe}
\t(\be+\ga)=\t(\be),\qqq \forall\, (\be,\ga)\in\cA \ts\G.
\]
This periodicity of the indices allows us to assign the \emph{index} $\t(\be_*)\in\Z^d$ to each oriented edge $\be_*\in\cA_*$ of the fundamental graph $\cG_*$ by setting
\[\lb{dco}
\t(\be_*)=\t(\be),
\]
where $\be\in\cA$ is an oriented edge in the periodic graph $\cG$ from the equivalence class $\be_*\in\cA_*=\cA/\G$. In other words, edge indices of the fundamental graph are induced by edge indices of the periodic graph. Due to (\ref{Gpe}), the edge index $\t(\be_*)$ is uniquely determined by \er{dco} and does not depend on the choice of $\be\in\cA$. For simplicity we use the same symbol $\t$ for the edge indices of both the periodic graph $\cG$ and its fundamental graph $\cG_*$. From the definition of the edge index it also follows that
\[\lb{inin}
\t(\ol\be\,)=-\t(\be), \qqq \forall\,\be\in\cA_*,
\]
where $\ol\be$ is the inverse edge of $\be$.

\begin{example}
The $d$-dimensional lattice $\mathbb{L}^d=(\mathcal{V},\mathcal{E})$, where the vertex set $\mathcal{V}$ and the edge set $\mathcal{E}$ are given by
$$
\mathcal{V}=\mathbb{Z}^d,\qqq
\mathcal{E}=\big\{(v,v+\mathfrak{a}_s), \;
\forall\,(v,s)\in\mathbb{Z}^d\times\mathbb{N}_d\big\},\qqq \mathbb{N}_d=\{1,\ldots,d\},
$$
and $\mathfrak{a}_1,\ldots,\mathfrak{a}_d$ is the standard basis of $\mathbb{Z}^d$ (for $d=2,3$ see Fig.~\ref{fig2}\emph{a},\emph{c}). The lattice $\mathbb{L}^d$ is a $\mathbb{Z}^d$-periodic graph with periods $\mathfrak{a}_1,\ldots,\mathfrak{a}_d$. The fundamental cell $\Omega=[0,1)^d$. The fundamental graph $\mathbb{L}^d_*$ of the lattice $\mathbb{L}^d$ consists of one vertex $v$ with degree $\vk_v=2d$ and $d$ loop edges $\mathbf{e}_1,\ldots,\mathbf{e}_d$ at this vertex $v$ with indices
$$
\tau(\mathbf{e}_1)=\mathfrak{a}_1, \quad \ldots\,,\quad \tau(\mathbf{e}_d)=\mathfrak{a}_d.
$$
\end{example}

\begin{example} The hexagonal lattice $\cG=(\cV,\cE)$ is shown in Fig.~\ref{figS1}\emph{a}, the vectors $\ga_1,\ga_2$ are the periods of $\cG$. The fundamental cell $\Omega$ is shaded in the figure. The vertex set $\cV$ and the edge set $\cE$ over the basis $\ga_1,\ga_2$ are given by
$$
\textstyle \cV=\Z^2\cup\big(\Z^2+\big(\frac13\,,\frac13\big)\big),
$$
$$
\textstyle\cE=\big\{\big(v,v+\big(\frac13\,,\frac13\big)\big),
\big(v,v+\big(-\frac23\,,\frac13\big)\big),\big(v,v+\big(\frac13\,,-\frac23\big)\big),
\quad\forall\,v\in\Z^2\big\}.
$$
The fundamental graph $\cG_*$ of $\cG$ consists of two vertices $v_1,v_2$ and three edges $\be_1,\be_2,\be_3$ connecting these vertices (Fig.~\ref{figS1}\emph{b}). The indices of the edges are given by
$$
\t(\be_1)=(0,0), \qq \t(\be_2)=(1,0), \qq \t(\be_3)=(0,1).
$$
\end{example}

\subsection{Discrete Schr\"odinger operators on periodic graphs.}
Let $\ell^2(\cV)$ be the Hilbert space of all square summable
functions  $f:\cV\to \C$ equipped with the norm
$$
\|f\|^2_{\ell^2(\cV)}=\sum_{v\in\cV}|f_v|^2<\infty.
$$

We consider the Schr\"odinger operator $H$ acting on the Hilbert
space $\ell^2(\cV)$ and given by
\[\lb{Sh}
H=\D+Q,
\]
where $\D$ is the discrete Laplacian having the form
\[\lb{ALO}
(\D f)_v=\sum_{(v,u)\in\cA}(f_v-f_u), \qqq f\in\ell^2(\cV), \qqq v\in\cV,
\]
and $Q$ is a real $\G$-periodic potential, i.e., it satisfies
$$
(Qf)_v=Q_vf_v, \qqq Q_{v+\ga}=Q_v, \qqq \forall\,(v,\ga)\in\cV\ts\G.
$$
The sum in (\ref{ALO}) is taken over all edges from $\mathcal{A}$ starting at
the vertex $v$. It is known that the Schr\"odinger operator $H$ is a bounded self-adjoint operator on $\ell^2(\cV)$ (see, e.g., \cite{SS92}) and the spectrum $\s(\D)$ of the Laplacian $\D$ satisfies \cite{MW89}
$$
0\in\s(\D)\subseteq[0,2\vk_+],\qqq\where \qqq
\vk_+=\max\limits_{v\in\cV_*}\vk_v.
$$

We briefly describe the spectrum of the Schr\"odinger operator $H$ on periodic graphs (for more details, see, e.g., \cite{HN09} or \cite{KS14}). We introduce the Hilbert space
$$
\mH=L^2\Big(\T^{d},{dk\/(2\pi)^d}\,;\ell^2(\cV_*)\Big)
=\int_{\T^{d}}^{\os}\ell^2(\cV_*)\,{dk \/(2\pi)^d}\,, \qqq
\T^d=\R^d/2\pi\Z^d=(-\pi,\pi]^d,
$$
equipped with the norm
$$
\|g\|^2_{\mH}=\int_{\T^d}\|g(k)\|_{\ell^2(\cV_*)}^2\frac{dk}{(2\pi)^d}\,,\qqq g\in\mH.
$$

The Schr\"odinger operator $H$ on $\ell^2(\cV)$ has the following decomposition into a constant fiber direct integral, see Theorem 1.1 from \cite{KS14},
$$
UH U^{-1}=\int^\oplus_{\T^d}H(k){dk\/(2\pi)^d}\,,
$$
where $U:\ell^2(\cV)\to\mH$ is some unitary operator (the Gelfand
transform), the parameter $k$ is called the \emph{quasimomentum}.
For each $k\in\T^d$ the fiber Schr\"odinger operator $H(k)$ on $\ell^2(\cV_*)$ is given by
\[\label{Hvt'}
H(k)=\D(k)+Q.
\]
Here $Q$ is the potential on $\ell^2(\cV_*)$, and $\D(k)$ is the fiber Laplacian having the form
\[
\label{fado}
\big(\D(k)f\big)_v=\sum_{\be=(v,u)\in\cA_*}\big(f_v-e^{i\lan\t(\be),\,k\ran}f_u\big),
 \qqq f\in\ell^2(\cV_*),\qqq v\in \cV_*,
\]
where $\t(\be)$ is the index of the edge $\be\in\cA_*$ defined by \er{in}, \er{dco}, and $\lan\cdot,\cdot\ran$ denotes the standard inner product in $\R^d$.

\begin{remark}\lb{Rfo0}
The fiber operator $H(0)$ is just the Schr\"odinger operator defined by \er{Sh}, \er{ALO} on the fundamental graph $\cG_*$.
\end{remark}

Let $\#M$ denote the number of elements in a set $M$. Each fiber operator $H(k)$, $k\in\T^{d}$, acts on the space $\ell^2(\cV_*)=\C^\n$, $\n=\#\cV_*$, and has $\n$ real eigenvalues $\l_j(k)$, $j=1,\ldots,\n$, labeled in non-decreasing order by
$$
\l_{1}(k)\leq\l_{2}(k)\leq\ldots\leq\l_{\nu}(k), \qqq
\forall\,k\in\T^{d},
$$
counting multiplicities.  Each \emph{band function} $\l_j(\cdot)$ is a continuous
and piecewise real analytic function on the torus $\T^{d}$ and creates the
\emph{spectral band} $\s_j(H)$ given by
\[\lb{ban.1H}
\begin{array}{l}
\s_j(H)=[\l_j^-,\l_j^+]=\l_j(\T^{d}),\qqq j\in\N_\n, \qqq \N_\n=\{1,\ldots,\n\},\\[6pt]
\displaystyle\textrm{where}\qqq \l_j^-=\min_{k\in\T^d}\l_j(k),\qqq \l_j^+=\max_{k\in\T^d}\l_j(k).
\end{array}
\]
Some of $\l_j(\cdot)$ may be constant, i.e., $\l_j(\cdot)=\L_j=\const$, on some subset of $\T^d$ of positive Lebesgue measure. In this case the Schr\"odinger operator $H$ on $\cG$ has the eigenvalue $\L_j$ of infinite multiplicity. We
call $\{\L_j\}$ a \emph{flat band}. Thus, the spectrum of the Schr\"odinger operator $H$ on a periodic graph $\cG$ has the form
\[\lb{specH}
\s(H)=\bigcup_{k\in\T^d}\s\big(H(k)\big)=
\bigcup_{j=1}^{\nu}\s_j(H)=\s_{ac}(H)\cup \s_{fb}(H),
\]
where $\s_{ac}(H)$ is the absolutely continuous spectrum, which is a
union of non-degenerate bands from \er{ban.1H}, and $\s_{fb}(H)$ is
the set of all flat bands.

\begin{remark}\lb{prbf}
\emph{i}) It is known (see, e.g., Theorem 2.1.\emph{ii} in \cite{KS19}) that the first spectral band $\s_1(H)=[\l_1^-,\l_1^+]$ is never a flat band and its lower point $\l_1^{-}=\l_1(0)$.

\emph{ii}) The band functions $\l_j(\cdot)$, $j\in\N_\n$, are even, i.e.,
$$
\l_j(-k)=\l_j(k),\qqq \forall\,k\in\T^d,
$$
see, e.g., Lemma 2.9.\emph{i} in \cite{KKR17}.
\end{remark}

\subsection{Perturbations of periodic graphs by adding edges}
Let $\G$ be a lattice of rank $d$ in $\R^d$ with basis $\ga_1,\ldots,\ga_d$, and let $\cG$ be a $\G$-periodic graph with the fundamental graph $\cG_*=\cG/\G$.

Let $t\in\Z^d$. We consider the $\G$-periodic graph $\cG_t$ obtained from the graph $\cG$ by adding an edge $\be_o$ with index $\t(\be_o)=t$ to its fundamental graph $\cG_*=(\cV_*,\cE_*)$ without changing the vertex set $\cV_*$. Recall that each edge of the fundamental graph $\cG_*$ is an equivalence class of edges of the periodic graph $\cG$. Then the edge $\be_o$ added to the fundamental graph $\cG_*$ gives an infinite number of $\G$-equivalent edges added to the periodic graph $\cG$ with the same index $\t(\be_o)$.

The $\G$-periodic graphs $\cG$ and $\cG_t$ will be called the \emph{unperturbed} and \emph{perturbed} graphs, respectively.

\begin{remark}
Each unoriented edge of the fundamental graph $\cG_*$ corresponds to two oppositely directed edges of $\cG_*$. Thus, the edge $\be_o$ added to the fundamental graph $\cG_*$ is considered as two oriented edges $\be_o,\ol\be_o$ with indices $\t(\be_o)=-\t(\ol\be_o)$, see (\ref{inin}).
\end{remark}

For example, as an unperturbed graph, let us consider the square lattice $\mathbb{L}^2$ with the minimal fundamental graph $\dL^2_*$ (Fig.~\ref{fig2}\emph{a}). Let $t\in\Z^2$. Then the perturbed graph $\dL^2_t$ is obtained from $\dL^2$ by adding an edge $(v,v+t)$ with index $t$ at each vertex $v$ of $\dL^2$. The perturbed graph $\dL^2_t$ when $t=(3,1)$ and its fundamental graph $\wt\dL^2_*$ are shown in  Fig.~\ref{fig2}\emph{b}. The added edges are marked in red color.

\smallskip

\begin{figure}[t!]\centering
\unitlength 1.1mm 
\linethickness{0.4pt}
\hspace{-8mm}
\begin{picture}(42,45)(0,0)
\hspace{-5mm}
\put(7,10){\line(1,0){36.00}}
\put(7,20){\line(1,0){36.00}}
\put(7,30){\line(1,0){36.00}}
\put(7,40){\line(1,0){36.00}}
\put(10,7){\line(0,1){36.00}}
\put(20,7){\line(0,1){36.00}}
\put(30,7){\line(0,1){36.00}}
\put(40,7){\line(0,1){36.00}}

\put(10,10){\vector(1,0){10.00}}
\put(10,10){\vector(0,1){10.00}}
\put(7.5,7.5){\small$v$}
\put(14,7.7){\small$\ga_1$}
\put(6.2,17){\small$\ga_2$}
\put(14,11){\small$\be_1$}
\put(10.5,15){\small$\be_2$}
\put(16,16){$\Omega$}
\put(5,42){$\dL^2$}

\bezier{15}(10.5,10)(10.5,15)(10.5,20)
\bezier{15}(11,10)(11,15)(11,20)
\bezier{15}(11.5,10)(11.5,15)(11.5,20)
\bezier{15}(12,10)(12,15)(12,20)
\bezier{15}(12.5,10)(12.5,15)(12.5,20)
\bezier{15}(13,10)(13,15)(13,20)
\bezier{15}(13.5,10)(13.5,15)(13.5,20)
\bezier{15}(14,10)(14,15)(14,20)
\bezier{15}(14.5,10)(14.5,15)(14.5,20)
\bezier{15}(15,10)(15,15)(15,20)
\bezier{15}(15.5,10)(15.5,15)(15.5,20)
\bezier{15}(16,10)(16,15)(16,20)
\bezier{15}(16.5,10)(16.5,15)(16.5,20)
\bezier{15}(17,10)(17,15)(17,20)
\bezier{15}(17.5,10)(17.5,15)(17.5,20)
\bezier{15}(18,10)(18,15)(18,20)
\bezier{15}(18.5,10)(18.5,15)(18.5,20)
\bezier{15}(19,10)(19,15)(19,20)
\bezier{15}(19.5,10)(19.5,15)(19.5,20)

\put(10,10){\circle*{1.3}}
\put(20,10){\circle*{1.3}}
\put(30,10){\circle*{1.3}}
\put(40,10){\circle*{1.3}}

\put(10,20){\circle*{1.3}}
\put(20,20){\circle*{1.3}}
\put(30,20){\circle*{1.3}}
\put(40,20){\circle*{1.3}}

\put(10,30){\circle*{1.3}}
\put(20,30){\circle*{1.3}}
\put(30,30){\circle*{1.3}}
\put(40,30){\circle*{1.3}}

\put(10,40){\circle*{1.3}}
\put(20,40){\circle*{1.3}}
\put(30,40){\circle*{1.3}}
\put(40,40){\circle*{1.3}}
\end{picture}\hspace{1mm}
\begin{picture}(42,45)
\hspace{-5mm}
\put(14,10.9){\small$\be_1$}
\put(10.5,15){\small$\be_2$}
\put(4,42){$\dL^2_t$}
\put(7,10){\line(1,0){36.00}}
\put(7,20){\line(1,0){36.00}}
\put(7,30){\line(1,0){36.00}}
\put(7,40){\line(1,0){36.00}}
\put(10,7){\line(0,1){36.00}}
\put(20,7){\line(0,1){36.00}}
\put(30,7){\line(0,1){36.00}}
\put(40,7){\line(0,1){36.00}}

\put(10,10){\vector(1,0){10.00}}
\put(10,10){\vector(0,1){10.00}}
\put(14,7.7){\small$\ga_1$}
\put(6.2,17){\small$\ga_2$}
\put(15,16){$\Omega$}
\put(7.5,7.5){\small$v$}
\put(40.0,17.8){\footnotesize$v\!+\!t$}

\bezier{15}(10.5,10)(10.5,15)(10.5,20)
\bezier{15}(11,10)(11,15)(11,20)
\bezier{15}(11.5,10)(11.5,15)(11.5,20)
\bezier{15}(12,10)(12,15)(12,20)
\bezier{15}(12.5,10)(12.5,15)(12.5,20)
\bezier{15}(13,10)(13,15)(13,20)
\bezier{15}(13.5,10)(13.5,15)(13.5,20)
\bezier{15}(14,10)(14,15)(14,20)
\bezier{15}(14.5,10)(14.5,15)(14.5,20)
\bezier{15}(15,10)(15,15)(15,20)
\bezier{15}(15.5,10)(15.5,15)(15.5,20)
\bezier{15}(16,10)(16,15)(16,20)
\bezier{15}(16.5,10)(16.5,15)(16.5,20)
\bezier{15}(17,10)(17,15)(17,20)
\bezier{15}(17.5,10)(17.5,15)(17.5,20)
\bezier{15}(18,10)(18,15)(18,20)
\bezier{15}(18.5,10)(18.5,15)(18.5,20)
\bezier{15}(19,10)(19,15)(19,20)
\bezier{15}(19.5,10)(19.5,15)(19.5,20)

\put(10,10){\circle*{1.3}}
\put(20,10){\circle*{1.3}}
\put(30,10){\circle*{1.3}}
\put(40,10){\circle*{1.3}}

\put(10,20){\circle*{1.3}}
\put(20,20){\circle*{1.3}}
\put(30,20){\circle*{1.3}}
\put(40,20){\circle*{1.3}}

\put(10,30){\circle*{1.3}}
\put(20,30){\circle*{1.3}}
\put(30,30){\circle*{1.3}}
\put(40,30){\circle*{1.3}}

\put(10,40){\circle*{1.3}}
\put(20,40){\circle*{1.3}}
\put(30,40){\circle*{1.3}}
\put(40,40){\circle*{1.3}}
\color{red}
\put(33.5,16.1){\small$\be_o$}
\put(10,10){\vector(3,1){30.00}}
\put(7,9){\line(3,1){36.00}}
\put(11,7){\line(3,1){32.00}}
\put(21,7){\line(3,1){22.00}}
\put(31,7){\line(3,1){12.00}}
\put(41,7){\line(3,1){2.00}}
\put(7,12.33){\line(3,1){36.00}}
\put(7,15.66){\line(3,1){36.00}}
\put(7,19){\line(3,1){36.00}}
\put(7,22.33){\line(3,1){36.00}}
\put(7,25.66){\line(3,1){36.00}}
\put(7,29){\line(3,1){36.00}}
\put(7,32.33){\line(3,1){32.00}}
\put(7,35.66){\line(3,1){22.00}}
\put(7,39){\line(3,1){12.00}}
\put(7,39.33){\line(3,1){2.00}}
\end{picture}\hspace{1mm}
\begin{picture}(42,45)(0,0)
\hspace{-5mm}
\put(8,10){\line(1,0){34.00}}
\put(8,20){\line(1,0){34.00}}
\put(8,30){\line(1,0){34.00}}
\put(8,40){\line(1,0){34.00}}

\put(12,13){\line(1,0){34.00}}
\put(12,23){\line(1,0){34.00}}
\put(12,33){\line(1,0){34.00}}
\put(12,43){\line(1,0){34.00}}

\put(16,16){\line(1,0){34.00}}
\put(16,26){\line(1,0){34.00}}
\put(16,36){\line(1,0){34.00}}
\put(16,46){\line(1,0){34.00}}

\put(10,10){\vector(1,0){10.00}}
\put(10,10){\vector(4,3){4.00}}

\put(7.7,7.4){\small$v$}
\put(5.8,20.8){\footnotesize$v\!+\!\ga_3$}
\put(5,42){$\dL^3$}

\put(14,7.7){\small$\ga_1$}
\put(10.2,13.5){\footnotesize$\ga_2$}

\put(10,10){\circle*{1.3}}
\put(20,10){\circle*{1.3}}
\put(30,10){\circle*{1.3}}
\put(40,10){\circle*{1.3}}

\put(10,20){\circle*{1.3}}
\put(20,20){\circle*{1.3}}
\put(30,20){\circle*{1.3}}
\put(40,20){\circle*{1.3}}

\put(10,30){\circle*{1.3}}
\put(20,30){\circle*{1.3}}
\put(30,30){\circle*{1.3}}
\put(40,30){\circle*{1.3}}

\put(10,40){\circle*{1.3}}
\put(20,40){\circle*{1.3}}
\put(30,40){\circle*{1.3}}
\put(40,40){\circle*{1.3}}

\put(14,13){\circle*{1.3}}
\put(24,13){\circle*{1.3}}
\put(34,13){\circle*{1.3}}
\put(44,13){\circle*{1.3}}

\put(14,23){\circle*{1.3}}
\put(24,23){\circle*{1.3}}
\put(34,23){\circle*{1.3}}
\put(44,23){\circle*{1.3}}

\put(14,33){\circle*{1.3}}
\put(24,33){\circle*{1.3}}
\put(34,33){\circle*{1.3}}
\put(44,33){\circle*{1.3}}

\put(14,43){\circle*{1.3}}
\put(24,43){\circle*{1.3}}
\put(34,43){\circle*{1.3}}
\put(44,43){\circle*{1.3}}

\put(18,16){\circle*{1.3}}
\put(28,16){\circle*{1.3}}
\put(38,16){\circle*{1.3}}
\put(48,16){\circle*{1.3}}

\put(18,26){\circle*{1.3}}
\put(28,26){\circle*{1.3}}
\put(38,26){\circle*{1.3}}
\put(48,26){\circle*{1.3}}

\put(18,36){\circle*{1.3}}
\put(28,36){\circle*{1.3}}
\put(38,36){\circle*{1.3}}
\put(48,36){\circle*{1.3}}

\put(18,46){\circle*{1.3}}
\put(28,46){\circle*{1.3}}
\put(38,46){\circle*{1.3}}
\put(48,46){\circle*{1.3}}

\put(8.4,8.8){\line(4,3){11.2}}
\put(18.4,8.8){\line(4,3){11.2}}
\put(28.4,8.8){\line(4,3){11.2}}
\put(38.4,8.8){\line(4,3){11.2}}
\put(8.4,18.8){\line(4,3){11.2}}
\put(18.4,18.8){\line(4,3){11.2}}
\put(28.4,18.8){\line(4,3){11.2}}
\put(38.4,18.8){\line(4,3){11.2}}
\put(8.4,28.8){\line(4,3){11.2}}
\put(18.4,28.8){\line(4,3){11.2}}
\put(28.4,28.8){\line(4,3){11.2}}
\put(38.4,28.8){\line(4,3){11.2}}
\put(8.4,38.8){\line(4,3){11.2}}
\put(18.4,38.8){\line(4,3){11.2}}
\put(28.4,38.8){\line(4,3){11.2}}
\put(38.4,38.8){\line(4,3){11.2}}

\put(6.2,17){\small$\ga_3$}
\put(14.5,16.5){$\Omega$}

\bezier{15}(10.5,10)(10.5,15)(10.5,20)
\bezier{15}(11,10)(11,15)(11,20)
\bezier{15}(11.5,10)(11.5,15)(11.5,20)
\bezier{15}(12,10)(12,15)(12,20)
\bezier{15}(12.5,10)(12.5,15)(12.5,20)
\bezier{15}(13,10)(13,15)(13,20)
\bezier{15}(13.5,10)(13.5,15)(13.5,20)
\bezier{15}(14,10)(14,15)(14,20)
\bezier{15}(14.5,10)(14.5,15)(14.5,20)
\bezier{15}(15,10)(15,15)(15,20)
\bezier{15}(15.5,10)(15.5,15)(15.5,20)
\bezier{15}(16,10)(16,15)(16,20)
\bezier{15}(16.5,10)(16.5,15)(16.5,20)
\bezier{15}(17,10)(17,15)(17,20)
\bezier{15}(17.5,10)(17.5,15)(17.5,20)
\bezier{15}(18,10)(18,15)(18,20)
\bezier{15}(18.5,10)(18.5,15)(18.5,20)
\bezier{15}(19,10)(19,15)(19,20)
\bezier{15}(19.5,10)(19.5,15)(19.5,20)

\bezier{15}(20.4,10.3)(20.4,15.3)(20.4,20.3)
\bezier{15}(20.8,10.6)(20.8,15.6)(20.8,20.6)
\bezier{15}(21.2,10.9)(21.2,15.9)(21.2,20.9)
\bezier{15}(21.6,11.2)(21.6,16.2)(21.6,21.2)
\bezier{15}(22.0,11.5)(22.0,16.5)(22.0,21.5)
\bezier{15}(22.4,11.8)(22.4,16.8)(22.4,21.8)
\bezier{15}(22.8,12.1)(22.8,17.1)(22.8,22.1)
\bezier{15}(23.2,12.4)(23.2,17.4)(23.2,22.4)
\bezier{15}(23.6,12.7)(23.6,17.7)(23.6,22.7)

\bezier{15}(10.4,20.3)(15.4,20.3)(20.4,20.3)
\bezier{15}(10.8,20.6)(15.8,20.6)(20.8,20.6)
\bezier{15}(11.2,20.9)(16.2,20.9)(21.2,20.9)
\bezier{15}(11.6,21.2)(16.6,21.2)(21.6,21.2)
\bezier{15}(12.0,21.5)(17.0,21.5)(22.0,21.5)
\bezier{15}(12.4,21.8)(17.4,21.8)(22.4,21.8)
\bezier{15}(12.8,22.1)(17.8,22.1)(22.8,22.1)
\bezier{15}(13.2,22.4)(18.2,22.4)(23.2,22.4)
\bezier{15}(13.6,22.7)(18.6,22.7)(23.6,22.7)

\color{red}
\put(10,8){\line(0,1){34.00}}
\put(20,8){\line(0,1){34.00}}
\put(30,8){\line(0,1){34.00}}
\put(40,8){\line(0,1){34.00}}
\put(14,11){\line(0,1){34.00}}
\put(24,11){\line(0,1){34.00}}
\put(34,11){\line(0,1){34.00}}
\put(44,11){\line(0,1){34.00}}
\put(18,14){\line(0,1){34.00}}
\put(28,14){\line(0,1){34.00}}
\put(38,14){\line(0,1){34.00}}
\put(48,14){\line(0,1){34.00}}

\put(6.2,13.8){\small$\be_o$}
\put(10,10){\vector(0,1){10.00}}

\end{picture}

\begin{picture}(35,22)
\put(-3,2){\emph{a})}
\put(2,18){$\dL^2_*$}
\put(27,4.5){$\be_2$}
\put(5,4.5){$\be_1$}
\put(17,1.5){$v$}
\put(23,10){\footnotesize $(0,1)$}
\put(4.5,10){\footnotesize $(1,0)$}
\put(18.0,5){\circle*{1.3}}
\bezier{200}(18.0,5)(5,12)(4.0,5)
\bezier{200}(18.0,5)(5,-2)(4.0,5)
\bezier{200}(18.0,5)(31,12)(32.0,5)
\bezier{200}(18.0,5)(31,-1)(32.0,5)
\end{picture}\hspace{8mm}
\begin{picture}(30,22)
\put(-3,2){\emph{b})}
\put(2,18){$\wt\dL^2_*$}
\put(23,10){\footnotesize $(0,1)$}
\put(4.5,10){\footnotesize $(1,0)$}

\put(27,4.5){$\be_2$}
\put(5,4.5){$\be_1$}
\put(17,1.5){$v$}
\put(18.0,5){\circle*{1.3}}
\bezier{200}(18.0,5)(5,12)(4.0,5)
\bezier{200}(18.0,5)(5,-2)(4.0,5)
\bezier{200}(18.0,5)(31,12)(32.0,5)
\bezier{200}(18.0,5)(31,-1)(32.0,5)

\color{red}
\put(13.5,20.3){\footnotesize $(t_1,t_2)$}
\put(16.5,16){$\be_o$}
\bezier{200}(18.0,5)(25,18)(18.0,19)
\bezier{200}(18.0,5)(11,18)(18.0,19)
\end{picture}
\hspace{16mm}
\begin{picture}(30,22)
\put(-4,2){\emph{c})}
\put(6,18){$\dL^3_*$}
\put(27,4.5){$\be_2$}
\put(5,4.5){$\be_1$}
\put(17,1.5){$v$}
\put(23,10){\footnotesize $(0,1,0)$}
\put(3,10){\footnotesize $(1,0,0)$}
\put(18.0,5){\circle*{1.3}}
\bezier{200}(18.0,5)(5,12)(4.0,5)
\bezier{200}(18.0,5)(5,-2)(4.0,5)
\bezier{200}(18.0,5)(31,12)(32.0,5)
\bezier{200}(18.0,5)(31,-1)(32.0,5)
\color{red}
\put(16.5,16){$\be_o$}
\put(12.5,20){\footnotesize $(0,0,1)$}
\bezier{200}(18.0,5)(25,18)(18.0,19)
\bezier{200}(18.0,5)(11,18)(18.0,19)
\end{picture}
\caption{\scriptsize\emph{a}) The square lattice $\dL^2$; $\ga_1,\ga_2$ are the periods of $\dL^2$; the fundamental cell $\Omega=[0,1)^2$ is shaded; the fundamental graph $\dL^2_*$ of $\dL^2$; the indices are shown near the edges. \emph{b}) The perturbed square lattice $\dL_t^2$; the added edges are marked in red color; the fundamental graph $\wt\dL^2_*$ of $\dL^2_t$; the red added edge $\be_o$ has the index $t=(t_1,t_2)\in\Z^2$. \emph{c})~The cubic lattice $\dL^3$; $\ga_1,\ga_2,\ga_3$ are the periods of $\dL^3$; the fundamental graph $\dL^3_*$ of $\dL^3$. The cubic lattice is obtained by stacking together infinitely many black copies of the square lattice $\dL^2$ along the vector $\ga_3$. The copies are connected by the red vertical edges $(v,v+\ga_3)$, $\forall v\in\Z^3$.} \label{fig2}
\end{figure}
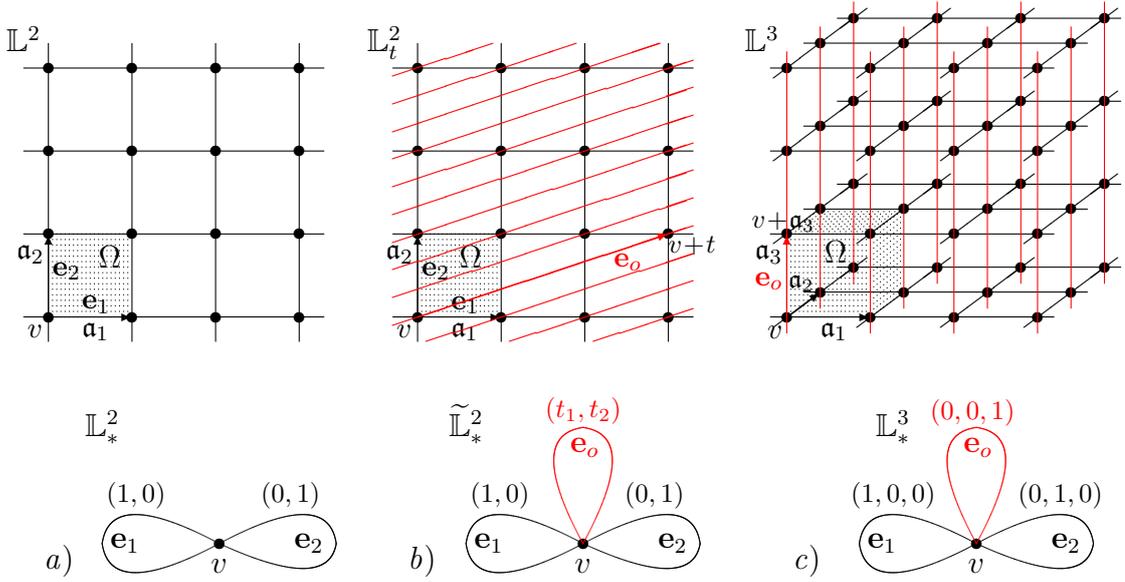

For each $t\in\Z^d$ the perturbed graph $\cG_t$ has the same
fundamental graph
$$
\wt\cG_*=\cG_t/\G=(\cV_*,\wt\cE_*),\qqq\textrm{where}\qqq
\wt\cE_*=\cE_*\cup\{\be_o\}.
$$

Due to (\ref{ban.1H}), (\ref{specH}), for each $t\in\Z^d$ the spectrum of the Schr\"odinger operator $H_{\cG_t}$ on the perturbed graph $\cG_t$ consists of $\nu$ bands $\s_j(H_{\cG_t})$:
$$
\s(H_{\cG_t})=\bigcup\limits_{j=1}^\nu\s_j(H_{\cG_t}), \qqq
\s_j(H_{\cG_t})=\big[\l_j^-(t),\l_j^+(t)\big], \qqq \nu=\#\cV_*.
$$
The \emph{band edges} $\l_j^\pm(t)$ depend on the parameter $t\in\Z^d$. Denote by $|\cdot|$ the Euclidean norm on $\R^d$.

We are interested in the following questions.

\smallskip

$\bullet$ Are there $\lim\limits_{|t|\to\infty}\lambda_j^\pm(t)$,\; $j\in\mathbb{N}_\nu$? And if so,

$\bullet$ Is there a periodic graph $\wt\cG$ such that the spectrum of the Schr\"odinger operator $H_{\wt\cG}=\D_{\wt\cG}+Q$ with a periodic potential $Q$ on $\wt\cG$ has the form
\[\lb{splg}
\begin{array}{l}
\displaystyle\s(H_{\wt\cG})=
\bigcup\limits_{j=1}^\nu\s_j(H_{\wt\cG}), \qqq
\s_j(H_{\wt\cG})=
\big[\,\wt\l_j^-,\wt\l_j^+\big],\\[12pt]
\textrm{where}\qqq \wt\l_j^-=\lim\limits_{|t|\to\iy}\l_j^-(t),\qqq
\wt\l_j^+=\lim\limits_{|t|\to\iy}\l_j^+(t)?
\end{array}
\]
If such a graph $\wt\cG$ exists, we will say that the perturbed graph $\cG_t$ is \emph{asymptotically isospectral} to the graph $\widetilde{\mathcal{G}}$ as $|t|\to\iy$.

\begin{remark}
The condition $|t|\to\iy$ means that the length of the edges added to the unperturbed periodic graph $\cG$ tends to infinity.
\end{remark}

The paper is organized as follows. In Section \ref{Sec2} we formulate our main results:

$\bullet$ we show that for any perturbed graph $\cG_t$ the graph $\wt\cG$ satisfying \er{splg} exists and give a construction of $\wt\cG$ (Theorem \ref{TAig});

$\bu$ we formulate a criterion for the perturbed graph $\cG_t$ to be not only asymptotically isospectral but just isospectral to the constructed graph $\wt\cG$ (Theorem \ref{TJig} and Corollary \ref{CJig});

$\bullet$ we obtain asymptotics of the band edges $\l_j^\pm(t)$, $j\in\N_\n$, of the Schr\"odinger operator $H_{\cG_t}$ on the perturbed graph $\cG_t$ as $|t|\to\iy$ (Theorem \ref{TAsy}).

\medskip

\no We illustrate the obtained results by some examples of asymptotically isospectral periodic graphs. In particular, we show that the 2-periodic Schr\"odinger operator on the one-dimensional lattice perturbed by long edges is asymptotically isospectral to the Schr\"odinger operator on the hexagonal lattice (Example \ref{Ehex}), and the $d$-dimensional lattice perturbed by long edges is asymptotically isospectral to the $(d+1)$-dimensional lattice (Example \ref{ExSl}).

In Section \ref{Sec3} we prove our main results. The proofs are based on the relation between the fiber Schr\"odinger operators for the perturbed graph $\cG_t$ and for the graph $\wt\cG$ (see Proposition \ref{Pnin0}). Section \ref{Sec4} is devoted to examples of asymptotically isospectral and just isospectral periodic graphs.

\section{Main results} \lb{Sec2}
\setcounter{equation}{0}

In order to answer the questions stated above, we proceed in the following way. First, for any perturbed graph $\cG_t$ we construct a graph $\wt\cG$ which will be called a \emph{limit} graph for $\cG_t$. Second, we show that the perturbed graph $\cG_t$ is asymptotically isospectral to its limit graph~$\wt\cG$.

\subsection{Limit graphs} \lb{Sec2.1} We provide a step-by-step construction of a limit graph $\wt\cG$. Each step will be  illustrated with the example of the hexagonal lattice as an unperturbed graph $\cG$.

\smallskip

\emph{Construction of the limit graph.}

\smallskip

\emph{\underline{Step 1.} Unperturbed graph $\cG$.} Let $\cG=(\cV,\cE)$ be a $\G$-periodic graph, where $\G$ is a lattice in $\R^d$ with basis $\ga_1,\ldots,\ga_d$ and the fundamental cell $\Omega$, and let $\cG_*=\cG/\G=(\cV_*,\cE_*)$ be its fundamental graph, see Fig.~\ref{figS1}. We identify the vertices of the fundamental graph $\cG_*$ with the corresponding vertices of the periodic graph $\cG$ from the fundamental cell $\Omega$.

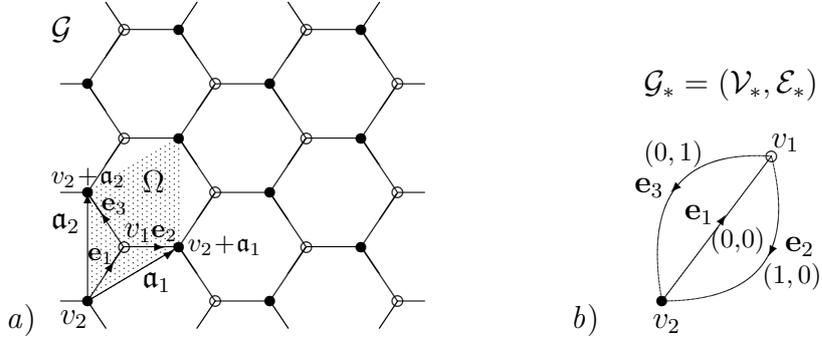
\begin{figure}[h!]\centering
\unitlength 1.20mm 
\linethickness{0.4pt}
\begin{picture}(35,35)(0,0)
\bezier{20}(4.5,3.3)(4.5,9.3)(4.5,15.3)
\bezier{20}(5.0,3.6)(5.0,9.6)(5.0,15.6)
\bezier{20}(5.5,3.9)(5.5,9.9)(5.5,15.9)
\bezier{20}(6.0,4.2)(6.0,10.2)(6.0,16.2)
\bezier{20}(6.5,4.5)(6.5,10.5)(6.5,16.5)
\bezier{20}(7.0,4.8)(7.0,10.8)(7.0,16.8)
\bezier{20}(7.5,5.1)(7.5,11.1)(7.5,17.1)
\bezier{20}(8.0,5.4)(8.0,11.4)(8.0,17.4)
\bezier{20}(8.5,5.7)(8.5,11.7)(8.5,17.7)
\bezier{20}(9.0,6.0)(9.0,12.0)(9.0,18.0)

\bezier{20}(9.5,6.3)(9.5,12.3)(9.5,18.3)
\bezier{20}(10.0,6.6)(10.0,12.6)(10.0,18.6)
\bezier{20}(10.5,6.9)(10.5,12.9)(10.5,18.9)
\bezier{20}(11.0,7.2)(11.0,13.2)(11.0,19.2)
\bezier{20}(11.5,7.5)(11.5,13.5)(11.5,19.5)
\bezier{20}(12.0,7.8)(12.0,13.8)(12.0,19.8)
\bezier{20}(12.5,8.1)(12.5,14.1)(12.5,20.1)
\bezier{20}(13.0,8.4)(13.0,14.4)(13.0,20.4)
\bezier{20}(13.5,8.7)(13.5,14.7)(13.5,20.7)
\bezier{20}(14,9)(14,15)(14,21)

\put(4,3){\circle*{1.3}}
\put(18,3){\circle{1.3}}
\put(24,3){\circle*{1.3}}
\put(38,3){\circle{1.3}}

\put(8,9){\circle{1.3}}
\put(14,9){\circle*{1.3}}
\put(28,9){\circle{1.3}}
\put(34,9){\circle*{1.3}}

\put(4,15){\circle*{1.3}}
\put(18,15){\circle{1.3}}
\put(24,15){\circle*{1.3}}
\put(38,15){\circle{1.3}}

\put(8,21){\circle{1.3}}
\put(14,21){\circle*{1.3}}
\put(28,21){\circle{1.3}}
\put(34,21){\circle*{1.3}}

\put(4,27){\circle*{1.3}}
\put(18,27){\circle{1.3}}
\put(24,27){\circle*{1.3}}
\put(38,27){\circle{1.3}}

\put(8,33){\circle{1.3}}
\put(14,33){\circle*{1.3}}
\put(28,33){\circle{1.3}}
\put(34,33){\circle*{1.3}}

\put(18,3){\line(1,0){6.00}}
\put(8,9){\line(1,0){6.00}}
\put(28,9){\line(1,0){6.00}}

\put(18,15){\line(1,0){6.00}}
\put(8,21){\line(1,0){6.00}}
\put(28,21){\line(1,0){6.00}}

\put(18,27){\line(1,0){6.00}}
\put(8,33){\line(1,0){6.00}}
\put(28,33){\line(1,0){6.00}}

\put(4,3){\line(2,3){4.00}}

\put(24,3){\line(2,3){4.00}}
\put(14,9){\line(2,3){4.00}}
\put(34,9){\line(2,3){4.00}}

\put(4,15){\line(2,3){4.00}}
\put(24,15){\line(2,3){4.00}}
\put(14,21){\line(2,3){4.00}}
\put(34,21){\line(2,3){4.00}}

\put(4,27){\line(2,3){4.00}}
\put(24,27){\line(2,3){4.00}}

\put(18,3){\line(-2,3){4.00}}
\put(38,3){\line(-2,3){4.00}}
\put(28,9){\line(-2,3){4.00}}
\put(8,9){\line(-2,3){4.00}}

\put(18,15){\line(-2,3){4.00}}
\put(38,15){\line(-2,3){4.00}}
\put(28,21){\line(-2,3){4.00}}
\put(8,21){\line(-2,3){4.00}}

\put(18,27){\line(-2,3){4.00}}
\put(38,27){\line(-2,3){4.00}}

\put(4,3){\vector(0,1){12.0}}
\put(13,8.3){\vector(3,2){0.5}}

\qbezier(4,3)(9,6)(14,9)
\put(4,3){\line(-1,0){3.00}}
\put(4,15){\line(-1,0){3.00}}
\put(4,27){\line(-1,0){3.00}}
\put(38,3){\line(1,0){3.00}}
\put(38,15){\line(1,0){3.00}}
\put(38,27){\line(1,0){3.00}}

\bezier{50}(4,3)(5,1.5)(6,0)
\bezier{50}(24,3)(25,1.5)(26,0)
\bezier{50}(18,3)(17,1.5)(16,0)
\bezier{50}(38,3)(37,1.5)(36,0)

\bezier{50}(8,33)(7,34.5)(6,36)
\bezier{50}(28,33)(27,34.5)(26,36)
\bezier{50}(14,33)(15,34.5)(16,36)
\bezier{50}(34,33)(35,34.5)(36,36)
\put(12,9){\vector(1,0){0.5}}
\put(6.4,6.6){\vector(2,3){0.5}}
\put(6,12){\vector(-2,3){0.5}}

\put(1,0.5){$v_2$}
\put(8,10){$v_1$}
\put(15,8.5){\footnotesize $v_2\!+\!\ga_1$}
\put(0,16){\footnotesize $v_2\!+\!\ga_2$}
\put(10,4.5){$\ga_1$}
\put(0.0,11){$\ga_2$}
\put(10,15){$\Omega$}
\put(0,32){$\cG$}
\put(11,10){\footnotesize  $\be_2$}
\put(5.5,13){\footnotesize $\be_3$}
\put(4.0,7.5){\footnotesize $\be_1$}
\put(-5,0){\emph{a})}
\end{picture}\hspace{25mm}
\begin{picture}(30,30)(0,0)
\put(0,0){\emph{b})}
\put(10,3){\circle*{1.3}}
\put(22,19){\circle{1.3}}
\put(9,0){$v_2$}
\put(22,20){$v_1$}

\bezier{200}(10,3.0)(16,11.0)(22,19.0)
\bezier{200}(10,3)(7,20)(22,19)
\bezier{200}(10,3)(26,3)(22,19)
\put(8,26){$\cG_*=(\cV_*,\cE_*)$}

\put(23.5,8.5){$\be_2$}
\put(7.0,15){$\be_3$}
\put(12.5,12){$\be_1$}

\put(16.0,11.0){\vector(2,3){1.0}}
\put(12.0,16.0){\vector(-1,-1){1.0}}
\put(22.7,10.0){\vector(-1,-2){1.0}}
\put(21,5){\footnotesize$(1,0)$}
\put(8,18.5){\footnotesize$(0,1)$}
\put(15.3,9){\footnotesize$(0,\!0)$}
\end{picture}
\caption{\scriptsize \emph{Step 1.} \emph{a}) An unperturbed graph $\cG$ (the hexagonal lattice as an example); $\ga_1,\ga_2$ are the periods of $\cG$; the fundamental cell $\Omega$ is shaded. \emph{b})~The fundamental graph $\cG_*$; the indices are shown near the edges.}
\label{figS1}
\end{figure}

\emph{\underline{Step 2.} Perturbed graph $\cG_t$.} Let $t=(t_1,\ldots,t_d)\in\Z^d$, and let $v_1,v_2\in\cV\cap\Omega$ be some fixed (maybe equal) vertices of the unperturbed graph $\cG$ from the fundamental cell $\Omega$. We consider the $\G$-periodic graph $\cG_t$ obtained from $\cG$ by adding infinitely many $\G$-equivalent edges
$$
\be_o+\ga,\qq \forall\,\ga\in\G,\qqq \textrm{where}\qqq
\be_o=(v_1,v_2+t_1\ga_1+\ldots+t_d\ga_d),
$$
Fig.~\ref{figS2}\emph{a}. The fundamental graph $\wt\cG_*=\cG_t/\G$ of the perturbed graph $\cG_t$ (Fig.~\ref{figS2}\emph{b}) is given by
\[\lb{pefg}
\wt\cG_*=(\cV_*,\wt\cE_*),\qqq\textrm{where}\qqq \wt\cE_*=\cE_*\cup\{\be_o\},\qq  \be_o=(v_1,v_2), \qq v_1,v_2\in\cV_*.
\]

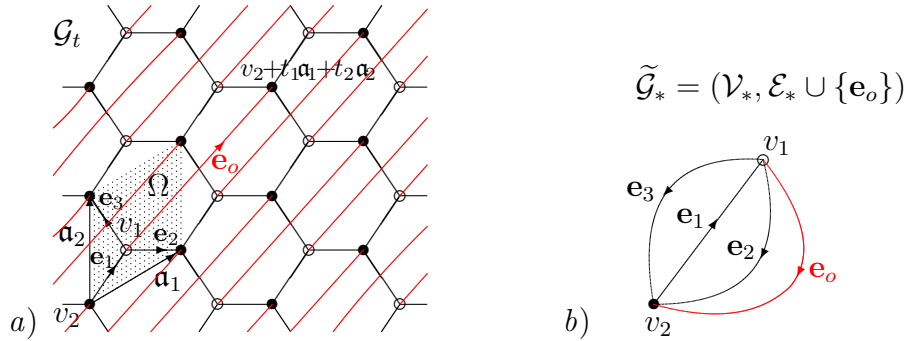
\begin{figure}[h!]\centering
\unitlength 1.20mm 
\linethickness{0.4pt}
\begin{picture}(35,35)(0,0)
\bezier{20}(4.5,3.3)(4.5,9.3)(4.5,15.3)
\bezier{20}(5.0,3.6)(5.0,9.6)(5.0,15.6)
\bezier{20}(5.5,3.9)(5.5,9.9)(5.5,15.9)
\bezier{20}(6.0,4.2)(6.0,10.2)(6.0,16.2)
\bezier{20}(6.5,4.5)(6.5,10.5)(6.5,16.5)
\bezier{20}(7.0,4.8)(7.0,10.8)(7.0,16.8)
\bezier{20}(7.5,5.1)(7.5,11.1)(7.5,17.1)
\bezier{20}(8.0,5.4)(8.0,11.4)(8.0,17.4)
\bezier{20}(8.5,5.7)(8.5,11.7)(8.5,17.7)
\bezier{20}(9.0,6.0)(9.0,12.0)(9.0,18.0)

\bezier{20}(9.5,6.3)(9.5,12.3)(9.5,18.3)
\bezier{20}(10.0,6.6)(10.0,12.6)(10.0,18.6)
\bezier{20}(10.5,6.9)(10.5,12.9)(10.5,18.9)
\bezier{20}(11.0,7.2)(11.0,13.2)(11.0,19.2)
\bezier{20}(11.5,7.5)(11.5,13.5)(11.5,19.5)
\bezier{20}(12.0,7.8)(12.0,13.8)(12.0,19.8)
\bezier{20}(12.5,8.1)(12.5,14.1)(12.5,20.1)
\bezier{20}(13.0,8.4)(13.0,14.4)(13.0,20.4)
\bezier{20}(13.5,8.7)(13.5,14.7)(13.5,20.7)
\bezier{20}(14,9)(14,15)(14,21)

\put(4,3){\circle*{1.3}}
\put(18,3){\circle{1.3}}
\put(24,3){\circle*{1.3}}
\put(38,3){\circle{1.3}}

\put(8,9){\circle{1.3}}
\put(14,9){\circle*{1.3}}
\put(28,9){\circle{1.3}}
\put(34,9){\circle*{1.3}}

\put(4,15){\circle*{1.3}}
\put(18,15){\circle{1.3}}
\put(24,15){\circle*{1.3}}
\put(38,15){\circle{1.3}}

\put(8,21){\circle{1.3}}
\put(14,21){\circle*{1.3}}
\put(28,21){\circle{1.3}}
\put(34,21){\circle*{1.3}}

\put(4,27){\circle*{1.3}}
\put(18,27){\circle{1.3}}
\put(24,27){\circle*{1.3}}
\put(38,27){\circle{1.3}}

\put(8,33){\circle{1.3}}
\put(14,33){\circle*{1.3}}
\put(28,33){\circle{1.3}}
\put(34,33){\circle*{1.3}}

\put(18,3){\line(1,0){6.00}}
\put(8,9){\line(1,0){6.00}}
\put(28,9){\line(1,0){6.00}}

\put(18,15){\line(1,0){6.00}}
\put(8,21){\line(1,0){6.00}}
\put(28,21){\line(1,0){6.00}}

\put(18,27){\line(1,0){6.00}}
\put(8,33){\line(1,0){6.00}}
\put(28,33){\line(1,0){6.00}}

\put(4,3){\line(2,3){4.00}}

\put(24,3){\line(2,3){4.00}}
\put(14,9){\line(2,3){4.00}}
\put(34,9){\line(2,3){4.00}}

\put(4,15){\line(2,3){4.00}}
\put(24,15){\line(2,3){4.00}}
\put(14,21){\line(2,3){4.00}}
\put(34,21){\line(2,3){4.00}}

\put(4,27){\line(2,3){4.00}}
\put(24,27){\line(2,3){4.00}}

\put(18,3){\line(-2,3){4.00}}
\put(38,3){\line(-2,3){4.00}}
\put(28,9){\line(-2,3){4.00}}
\put(8,9){\line(-2,3){4.00}}

\put(18,15){\line(-2,3){4.00}}
\put(38,15){\line(-2,3){4.00}}
\put(28,21){\line(-2,3){4.00}}
\put(8,21){\line(-2,3){4.00}}

\put(18,27){\line(-2,3){4.00}}
\put(38,27){\line(-2,3){4.00}}

\put(4,3){\vector(0,1){12.0}}
\put(13,8.3){\vector(3,2){0.5}}

\qbezier(4,3)(9,6)(14,9)

\put(11,5){$\ga_1$}
\put(0.5,10){$\ga_2$}
\put(0,1){$v_2$}
\put(7,11){$v_1$}
\put(20.5,28.2){\footnotesize$v_2\!\!+\!\!t_1\!\ga_{\!1}\!\!+\!\!t_2\ga_{\!2}$}
\put(10.3,15){$\Omega$}
\put(0,32){$\cG_t$}

\put(4,3){\line(-1,0){3.00}}
\put(4,15){\line(-1,0){3.00}}
\put(4,27){\line(-1,0){3.00}}
\put(38,3){\line(1,0){3.00}}
\put(38,15){\line(1,0){3.00}}
\put(38,27){\line(1,0){3.00}}

\bezier{50}(4,3)(5,1.5)(6,0)
\bezier{50}(24,3)(25,1.5)(26,0)
\bezier{50}(18,3)(17,1.5)(16,0)
\bezier{50}(38,3)(37,1.5)(36,0)

\bezier{50}(8,33)(7,34.5)(6,36)
\bezier{50}(28,33)(27,34.5)(26,36)
\bezier{50}(14,33)(15,34.5)(16,36)
\bezier{50}(34,33)(35,34.5)(36,36)
\put(12,9){\vector(1,0){0.5}}
\put(6.4,6.6){\vector(2,3){0.5}}
\put(6,12){\vector(-2,3){0.5}}

\put(11,10){\footnotesize$\be_2$}
\put(5.0,14){\footnotesize$\be_3$}
\put(4.0,7.5){\footnotesize$\be_1$}
\put(-5,0){\emph{a})}
\color{red}
\put(17.3,18){$\be_o$}
\put(17.8,20){\vector(1,1){1}}
\qbezier(4,3)(2.65,1.5)(1.3,0)
\qbezier(24,3)(22.65,1.5)(21.3,0)
\qbezier(14,9)(10,4.5)(6,0)
\qbezier(34,9)(30,4.5)(26,0)
\qbezier(4,15)(2,12.25)(0,10.5)
\qbezier(24,15)(17.35,7.5)(10.7,0)
\qbezier(14,21)(7,13.1)(0,5.25)
\qbezier(4,27)(2,24.25)(0,22.5)
\qbezier(14,33)(7,26.1)(0,17.25)
\qbezier(41,11.6)(35.15,5.8)(30.7,0)

\qbezier(18,3)(26,12)(34,21)
\qbezier(38,3)(39.5,4.7)(41,6.4)
\qbezier(8,9)(16,18)(24,27)
\qbezier(28,9)(34.5,16.3)(41,23.6)
\qbezier(18,15)(26,24)(34,33)
\qbezier(38,15)(39.5,16.7)(41,18.4)
\qbezier(8,21)(14.65,28.5)(21.3,36)
\qbezier(28,21)(34.5,28.3)(41,35.6)
\qbezier(18,27)(22,31.5)(26,36)
\qbezier(38,27)(39.5,28.7)(41,30.4)
\qbezier(8,33)(9.35,34.5)(10.7,36)
\qbezier(28,33)(29.35,34.5)(30.7,36)

\end{picture}\hspace{25mm}\begin{picture}(30,30)(0,0)
\put(0,0){\emph{b})}
\put(10,3){\circle*{1.3}}
\put(22,19){\circle{1.3}}
\put(9,0){$v_2$}
\put(22,20){$v_1$}

\bezier{200}(10,3.0)(16,11.0)(22,19.0)
\bezier{200}(10,3)(7,20)(22,19)
\bezier{200}(10,3)(26,3)(22,19)
\put(8,26){$\wt\cG_*=(\cV_*,\cE_*\cup\{\be_o\})$}
\put(18,8.5){$\be_2$}
\put(7.0,15){$\be_3$}
\put(12.5,12){$\be_1$}

\put(16.0,11.0){\vector(2,3){1.0}}
\put(12.0,16.0){\vector(-1,-1){1.0}}
\put(22.7,10.0){\vector(-1,-2){1.0}}
\color{red}
\put(27,6){$\be_o$}
\bezier{200}(10,3.0)(20,0.15)(25,5)
\bezier{200}(22,19.0)(29,10)(25,5)
\put(26.7,8.0){\vector(-1,-2){1.0}}
\end{picture}
\caption{\scriptsize \emph{Step 2.} \emph{a}) The perturbed graph $\cG_t$; the added edges are marked in red color. \emph{b}) The fundamental graph $\wt\cG_*$ of $\cG_t$; the red added edge $\be_o$ has the index $t=(t_1,t_2)\in\Z^2$.}
\label{figS2}
\end{figure}

\emph{\underline{Step 3.} Limit graph $\wt\cG$.} We complete the basis $\ga_1,\ldots,\ga_d$ of $\R^d$ to a basis of $\R^{d+1}$ by some vector $\ga_{d+1}$, and denote by $\wt\G$ the lattice in $\R^{d+1}$ with basis $\ga_1,\ldots,\ga_d,\ga_{d+1}$. The \emph{limit} graph $\wt\cG$ is obtained by the following way.

$\bu$ First, we stack together infinitely many copies
$$
\cG+n\ga_{d+1}\ss\R^{d+1},\qqq \forall\,n\in\Z,
$$
of the unperturbed graph $\cG\ss\R^d$ along the vector $\ga_{d+1}$, see Fig.~\ref{figS3}\emph{a}. Here $\R^d$ is considered as the subspace $\R^d\ts\{0\}$ of $\R^{d+1}$.

$\bu$ Second, these copies of $\cG$ are connected in a $\wt\G$-periodic way by the edges
$$
\be_o+\ga,\qq \forall\,\ga\in\wt\G,\qqq \textrm{where}\qqq
\be_o=(v_1,v_2+\ga_{d+1}),
$$
see Fig.~\ref{figS3}\emph{b}. The edge $\be_o$ connects the vertex $v_1$ of the "zero"\, copy $\cG$ with the vertex $v_2+\ga_{d+1}$ of the next copy $\cG+\ga_{d+1}$. Recall that $v_1,v_2$ are some fixed vertices of the unperturbed graph $\cG$ from its fundamental cell $\Omega$.

By construction, the obtained graph $\wt\cG$ is a $\wt\G$-periodic graph with the same fundamental graph $\wt\cG_*=\wt\cG/\wt\G$ as the perturbed graph $\cG_t$, see \er{pefg} and  Fig.~\ref{figS2}\emph{b}.

\begin{figure}[h]
\unitlength 1.0mm 
\linethickness{0.4pt}
\begin{picture}(70,60)
\hspace{-7mm}
\put(12,0){\vector(2,1){12.00}}
\put(12,0){\vector(1,0){18.00}}
\put(12,0){\vector(0,1){20.00}}
\bezier{30}(13.0,0.5)(22.0,0.5)(31.0,0.5)
\bezier{30}(14.0,1.0)(23.0,1.0)(32.0,1.0)
\bezier{30}(15.0,1.5)(24.0,1.5)(33.0,1.5)
\bezier{30}(16.0,2.0)(25.0,2.0)(34.0,2.0)
\bezier{30}(17.0,2.5)(26.0,2.5)(35.0,2.5)
\bezier{30}(18.0,3.0)(27.0,3.0)(36.0,3.0)
\bezier{30}(19.0,3.5)(28.0,3.5)(37.0,3.5)
\bezier{30}(20.0,4.0)(29.0,4.0)(38.0,4.0)
\bezier{30}(21.0,4.5)(30.0,4.5)(39.0,4.5)
\bezier{30}(22.0,5.0)(31.0,5.0)(40.0,5.0)
\bezier{30}(23.0,5.5)(32.0,5.5)(41.0,5.5)
\bezier{30}(24.0,6.0)(33.0,6.0)(42.0,6.0)
\put(0,-2.5){\emph{a})}
\put(30,2){$\Omega$}
\put(23,-2.5){$\ga_1$}
\put(15.5,4.5){$\ga_2$}
\put(8,13){$\ga_3$}
\put(20,0.5){\small $v_1$}
\put(7,-1.5){\small $v_2$}
\put(7.5,22){\small$v_2\!\!+\!\!\ga_3$}
\put(60,10){$\cG$}
\put(57,29.5){$\cG\!+\!\ga_3$}
\put(56.5,49.5){$\cG\!+\!2\ga_3$}

\qbezier(12,0)(11.5,-0.5)(11,-1)
\qbezier(30,0)(29.5,-0.5)(29,-1)
\qbezier(48,0)(47.5,-0.5)(47,-1)
\qbezier(3,3)(2.25,2.75)(1.5,2.5)
\qbezier(57,3)(57.75,2.75)(58.5,2.5)
\qbezier(6,6)(5.25,6.25)(4.5,6.5)
\qbezier(69,9)(69.75,8.75)(70.5,8.5)
\qbezier(18,12)(17.25,12.25)(16.5,12.5)
\qbezier(72,12)(72.75,12.25)(73.5,12.5)
\qbezier(27,15)(27.5,15.5)(28,16)
\qbezier(45,15)(45.5,15.5)(46,16)
\qbezier(63,15)(63.5,15.5)(64,16)

\qbezier(12,20)(11.5,19.5)(11,19)
\qbezier(30,20)(29.5,19.5)(29,19)
\qbezier(48,20)(47.5,19.5)(47,19)
\qbezier(3,23)(2.25,22.75)(1.5,22.5)
\qbezier(57,23)(57.75,22.75)(58.5,22.5)
\qbezier(6,26)(5.25,26.25)(4.5,26.5)
\qbezier(69,29)(69.75,28.75)(70.5,28.5)
\qbezier(18,32)(17.25,32.25)(16.5,32.5)
\qbezier(72,32)(72.75,32.25)(73.5,32.5)
\qbezier(27,35)(27.5,35.5)(28,36)
\qbezier(45,35)(45.5,35.5)(46,36)
\qbezier(63,35)(63.5,35.5)(64,36)

\qbezier(12,40)(11.5,39.5)(11,39)
\qbezier(30,40)(29.5,39.5)(29,39)
\qbezier(48,40)(47.5,39.5)(47,39)
\qbezier(3,43)(2.25,42.75)(1.5,42.5)
\qbezier(57,43)(57.75,42.75)(58.5,42.5)
\qbezier(6,46)(5.25,46.25)(4.5,46.5)
\qbezier(69,49)(69.75,48.75)(70.5,48.5)
\qbezier(18,52)(17.25,52.25)(16.5,52.5)
\qbezier(72,52)(72.75,52.25)(73.5,52.5)
\qbezier(27,55)(27.5,55.5)(28,56)
\qbezier(45,55)(45.5,55.5)(46,56)
\qbezier(63,55)(63.5,55.5)(64,56)

\put(12,0){\circle*{1.5}}
\put(30,0){\circle*{1.5}}
\put(48,0){\circle*{1.5}}
\put(3,3){\circle{1.5}}
\put(21,3){\circle{1.5}}
\put(39,3){\circle{1.5}}
\put(57,3){\circle{1.5}}
\put(6,6){\circle*{1.5}}
\put(24,6){\circle*{1.5}}
\put(42,6){\circle*{1.5}}
\put(60,6){\circle*{1.5}}
\put(15,9){\circle{1.5}}
\put(33,9){\circle{1.5}}
\put(51,9){\circle{1.5}}
\put(69,9){\circle{1.5}}
\put(18,12){\circle*{1.5}}
\put(36,12){\circle*{1.5}}
\put(54,12){\circle*{1.5}}
\put(72,12){\circle*{1.5}}
\put(27,15){\circle{1.5}}
\put(45,15){\circle{1.5}}
\put(63,15){\circle{1.5}}

\put(12,0){\line(-3,1){9.00}}
\put(12,0){\line(3,1){9.00}}
\put(30,0){\line(-3,1){9.00}}
\put(30,0){\line(3,1){9.00}}
\put(48,0){\line(-3,1){9.00}}
\put(48,0){\line(3,1){9.00}}

\qbezier(6,6)(4.5,4.5)(3,3)
\put(6,6){\line(3,1){9.00}}
\qbezier(24,6)(22.5,4.5)(21,3)
\put(24,6){\line(3,1){9.00}}
\put(24,6){\line(-3,1){9.00}}
\qbezier(42,6)(40.5,4.5)(39,3)
\put(42,6){\line(3,1){9.00}}
\put(42,6){\line(-3,1){9.00}}
\qbezier(60,6)(58.5,4.5)(57,3)
\put(60,6){\line(3,1){9.00}}
\put(60,6){\line(-3,1){9.00}}
\qbezier(18,12)(16.5,10.5)(15,9)
\put(18,12){\line(3,1){9.00}}
\qbezier(36,12)(34.5,10.5)(33,9)
\put(36,12){\line(3,1){9.00}}
\put(36,12){\line(-3,1){9.00}}
\qbezier(54,12)(52.5,10.5)(51,9)
\put(54,12){\line(3,1){9.00}}
\put(54,12){\line(-3,1){9.00}}
\qbezier(72,12)(70.5,10.5)(69,9)
\put(72,12){\line(-3,1){9.00}}


\put(12,20){\circle*{1.5}}
\put(30,20){\circle*{1.5}}
\put(48,20){\circle*{1.5}}
\put(3,23){\circle{1.5}}
\put(21,23){\circle{1.5}}
\put(39,23){\circle{1.5}}
\put(57,23){\circle{1.5}}
\put(6,26){\circle*{1.5}}
\put(24,26){\circle*{1.5}}
\put(42,26){\circle*{1.5}}
\put(60,26){\circle*{1.5}}
\put(15,29){\circle{1.5}}
\put(33,29){\circle{1.5}}
\put(51,29){\circle{1.5}}
\put(69,29){\circle{1.5}}
\put(18,32){\circle*{1.5}}
\put(36,32){\circle*{1.5}}
\put(54,32){\circle*{1.5}}
\put(72,32){\circle*{1.5}}
\put(27,35){\circle{1.5}}
\put(45,35){\circle{1.5}}
\put(63,35){\circle{1.5}}

\put(12,20){\line(-3,1){9.00}}
\put(12,20){\line(3,1){9.00}}
\put(30,20){\line(-3,1){9.00}}
\put(30,20){\line(3,1){9.00}}
\put(48,20){\line(-3,1){9.00}}
\put(48,20){\line(3,1){9.00}}

\qbezier(6,26)(4.5,24.5)(3,23)
\put(6,26){\line(3,1){9.00}}
\qbezier(24,26)(22.5,24.5)(21,23)
\put(24,26){\line(3,1){9.00}}
\put(24,26){\line(-3,1){9.00}}
\qbezier(42,26)(40.5,24.5)(39,23)
\put(42,26){\line(3,1){9.00}}
\put(42,26){\line(-3,1){9.00}}
\qbezier(60,26)(58.5,24.5)(57,23)
\put(60,26){\line(3,1){9.00}}
\put(60,26){\line(-3,1){9.00}}
\qbezier(18,32)(16.5,30.5)(15,29)
\put(18,32){\line(3,1){9.00}}
\qbezier(36,32)(34.5,30.5)(33,29)
\put(36,32){\line(3,1){9.00}}
\put(36,32){\line(-3,1){9.00}}
\qbezier(54,32)(52.5,30.5)(51,29)
\put(54,32){\line(3,1){9.00}}
\put(54,32){\line(-3,1){9.00}}
\qbezier(72,32)(70.5,30.5)(69,29)
\put(72,32){\line(-3,1){9.00}}


\put(12,40){\circle*{1.5}}
\put(30,40){\circle*{1.5}}
\put(48,40){\circle*{1.5}}
\put(3,43){\circle{1.5}}
\put(21,43){\circle{1.5}}
\put(39,43){\circle{1.5}}
\put(57,43){\circle{1.5}}
\put(6,46){\circle*{1.5}}
\put(24,46){\circle*{1.5}}
\put(42,46){\circle*{1.5}}
\put(60,46){\circle*{1.5}}
\put(15,49){\circle{1.5}}
\put(33,49){\circle{1.5}}
\put(51,49){\circle{1.5}}
\put(69,49){\circle{1.5}}
\put(18,52){\circle*{1.5}}
\put(36,52){\circle*{1.5}}
\put(54,52){\circle*{1.5}}
\put(72,52){\circle*{1.5}}
\put(27,55){\circle{1.5}}
\put(45,55){\circle{1.5}}
\put(63,55){\circle{1.5}}

\put(12,40){\line(-3,1){9.00}}
\put(12,40){\line(3,1){9.00}}
\put(30,40){\line(-3,1){9.00}}
\put(30,40){\line(3,1){9.00}}
\put(48,40){\line(-3,1){9.00}}
\put(48,40){\line(3,1){9.00}}

\qbezier(6,46)(4.5,44.5)(3,43)
\put(6,46){\line(3,1){9.00}}
\qbezier(24,46)(22.5,44.5)(21,43)
\put(24,46){\line(3,1){9.00}}
\put(24,46){\line(-3,1){9.00}}
\qbezier(42,46)(40.5,44.5)(39,43)
\put(42,46){\line(3,1){9.00}}
\put(42,46){\line(-3,1){9.00}}
\qbezier(60,46)(58.5,44.5)(57,43)
\put(60,46){\line(3,1){9.00}}
\put(60,46){\line(-3,1){9.00}}

\qbezier(18,52)(16.5,50.5)(15,49)
\put(18,52){\line(3,1){9.00}}
\qbezier(36,52)(34.5,50.5)(33,49)
\put(36,52){\line(3,1){9.00}}
\put(36,52){\line(-3,1){9.00}}
\qbezier(54,52)(52.5,50.5)(51,49)
\put(54,52){\line(3,1){9.00}}
\put(54,52){\line(-3,1){9.00}}
\qbezier(72,52)(70.5,50.5)(69,49)
\put(72,52){\line(-3,1){9.00}}
\end{picture}\hspace{1mm}
\begin{picture}(65,60)
\put(12,0){\vector(2,1){12.00}}
\put(12,0){\vector(1,0){18.00}}
\put(12,0){\vector(0,1){20.00}}

\bezier{30}(13.0,0.5)(22.0,0.5)(31.0,0.5)
\bezier{30}(14.0,1.0)(23.0,1.0)(32.0,1.0)
\bezier{30}(15.0,1.5)(24.0,1.5)(33.0,1.5)
\bezier{30}(16.0,2.0)(25.0,2.0)(34.0,2.0)
\bezier{30}(17.0,2.5)(26.0,2.5)(35.0,2.5)
\bezier{30}(18.0,3.0)(27.0,3.0)(36.0,3.0)
\bezier{30}(19.0,3.5)(28.0,3.5)(37.0,3.5)
\bezier{30}(20.0,4.0)(29.0,4.0)(38.0,4.0)
\bezier{30}(21.0,4.5)(30.0,4.5)(39.0,4.5)
\bezier{30}(22.0,5.0)(31.0,5.0)(40.0,5.0)
\bezier{30}(23.0,5.5)(32.0,5.5)(41.0,5.5)
\bezier{30}(24.0,6.0)(33.0,6.0)(42.0,6.0)
\put(0,-2.5){\emph{b})}
\put(30,2){$\Omega$}
\put(23,-2.5){$\ga_1$}
\put(15.5,4.5){$\ga_2$}
\put(8,13){$\ga_3$}
\put(20,0.5){\small $v_1$}
\put(7,-1.5){\small $v_2$}
\put(7.5,22){\small$v_2\!\!+\!\!\ga_3$}
\put(60,10){$\cG$}
\put(57,29.5){$\cG\!+\!\ga_3$}
\put(55.5,49.5){$\cG\!+\!2\ga_3$}
\put(6,53){$\wt\cG$}

\qbezier(12,0)(11.5,-0.5)(11,-1)
\qbezier(30,0)(29.5,-0.5)(29,-1)
\qbezier(48,0)(47.5,-0.5)(47,-1)
\qbezier(3,3)(2.25,2.75)(1.5,2.5)
\qbezier(57,3)(57.75,2.75)(58.5,2.5)
\qbezier(6,6)(5.25,6.25)(4.5,6.5)
\qbezier(69,9)(69.75,8.75)(70.5,8.5)
\qbezier(18,12)(17.25,12.25)(16.5,12.5)
\qbezier(72,12)(72.75,12.25)(73.5,12.5)
\qbezier(27,15)(27.5,15.5)(28,16)
\qbezier(45,15)(45.5,15.5)(46,16)
\qbezier(63,15)(63.5,15.5)(64,16)

\qbezier(12,20)(11.5,19.5)(11,19)
\qbezier(30,20)(29.5,19.5)(29,19)
\qbezier(48,20)(47.5,19.5)(47,19)
\qbezier(3,23)(2.25,22.75)(1.5,22.5)
\qbezier(57,23)(57.75,22.75)(58.5,22.5)
\qbezier(6,26)(5.25,26.25)(4.5,26.5)
\qbezier(69,29)(69.75,28.75)(70.5,28.5)
\qbezier(18,32)(17.25,32.25)(16.5,32.5)
\qbezier(72,32)(72.75,32.25)(73.5,32.5)
\qbezier(27,35)(27.5,35.5)(28,36)
\qbezier(45,35)(45.5,35.5)(46,36)
\qbezier(63,35)(63.5,35.5)(64,36)

\qbezier(12,40)(11.5,39.5)(11,39)
\qbezier(30,40)(29.5,39.5)(29,39)
\qbezier(48,40)(47.5,39.5)(47,39)
\qbezier(3,43)(2.25,42.75)(1.5,42.5)
\qbezier(57,43)(57.75,42.75)(58.5,42.5)
\qbezier(6,46)(5.25,46.25)(4.5,46.5)
\qbezier(69,49)(69.75,48.75)(70.5,48.5)
\qbezier(18,52)(17.25,52.25)(16.5,52.5)
\qbezier(72,52)(72.75,52.25)(73.5,52.5)
\qbezier(27,55)(27.5,55.5)(28,56)
\qbezier(45,55)(45.5,55.5)(46,56)
\qbezier(63,55)(63.5,55.5)(64,56)

\put(12,0){\circle*{1.5}}
\put(30,0){\circle*{1.5}}
\put(48,0){\circle*{1.5}}
\put(3,3){\circle{1.5}}
\put(21,3){\circle{1.5}}
\put(39,3){\circle{1.5}}
\put(57,3){\circle{1.5}}
\put(6,6){\circle*{1.5}}
\put(24,6){\circle*{1.5}}
\put(42,6){\circle*{1.5}}
\put(60,6){\circle*{1.5}}
\put(15,9){\circle{1.5}}
\put(33,9){\circle{1.5}}
\put(51,9){\circle{1.5}}
\put(69,9){\circle{1.5}}
\put(18,12){\circle*{1.5}}
\put(36,12){\circle*{1.5}}
\put(54,12){\circle*{1.5}}
\put(72,12){\circle*{1.5}}
\put(27,15){\circle{1.5}}
\put(45,15){\circle{1.5}}
\put(63,15){\circle{1.5}}

\put(12,0){\line(-3,1){9.00}}
\put(12,0){\line(3,1){9.00}}
\put(30,0){\line(-3,1){9.00}}
\put(30,0){\line(3,1){9.00}}
\put(48,0){\line(-3,1){9.00}}
\put(48,0){\line(3,1){9.00}}

\qbezier(6,6)(4.5,4.5)(3,3)
\put(6,6){\line(3,1){9.00}}
\qbezier(24,6)(22.5,4.5)(21,3)
\put(24,6){\line(3,1){9.00}}
\put(24,6){\line(-3,1){9.00}}
\qbezier(42,6)(40.5,4.5)(39,3)
\put(42,6){\line(3,1){9.00}}
\put(42,6){\line(-3,1){9.00}}
\qbezier(60,6)(58.5,4.5)(57,3)
\put(60,6){\line(3,1){9.00}}
\put(60,6){\line(-3,1){9.00}}
\qbezier(18,12)(16.5,10.5)(15,9)
\put(18,12){\line(3,1){9.00}}
\qbezier(36,12)(34.5,10.5)(33,9)
\put(36,12){\line(3,1){9.00}}
\put(36,12){\line(-3,1){9.00}}
\qbezier(54,12)(52.5,10.5)(51,9)
\put(54,12){\line(3,1){9.00}}
\put(54,12){\line(-3,1){9.00}}
\qbezier(72,12)(70.5,10.5)(69,9)
\put(72,12){\line(-3,1){9.00}}


\put(12,20){\circle*{1.5}}
\put(30,20){\circle*{1.5}}
\put(48,20){\circle*{1.5}}
\put(3,23){\circle{1.5}}
\put(21,23){\circle{1.5}}
\put(39,23){\circle{1.5}}
\put(57,23){\circle{1.5}}
\put(6,26){\circle*{1.5}}
\put(24,26){\circle*{1.5}}
\put(42,26){\circle*{1.5}}
\put(60,26){\circle*{1.5}}
\put(15,29){\circle{1.5}}
\put(33,29){\circle{1.5}}
\put(51,29){\circle{1.5}}
\put(69,29){\circle{1.5}}
\put(18,32){\circle*{1.5}}
\put(36,32){\circle*{1.5}}
\put(54,32){\circle*{1.5}}
\put(72,32){\circle*{1.5}}
\put(27,35){\circle{1.5}}
\put(45,35){\circle{1.5}}
\put(63,35){\circle{1.5}}

\put(12,20){\line(-3,1){9.00}}
\put(12,20){\line(3,1){9.00}}
\put(30,20){\line(-3,1){9.00}}
\put(30,20){\line(3,1){9.00}}
\put(48,20){\line(-3,1){9.00}}
\put(48,20){\line(3,1){9.00}}

\qbezier(6,26)(4.5,24.5)(3,23)
\put(6,26){\line(3,1){9.00}}
\qbezier(24,26)(22.5,24.5)(21,23)
\put(24,26){\line(3,1){9.00}}
\put(24,26){\line(-3,1){9.00}}
\qbezier(42,26)(40.5,24.5)(39,23)
\put(42,26){\line(3,1){9.00}}
\put(42,26){\line(-3,1){9.00}}
\qbezier(60,26)(58.5,24.5)(57,23)
\put(60,26){\line(3,1){9.00}}
\put(60,26){\line(-3,1){9.00}}
\qbezier(18,32)(16.5,30.5)(15,29)
\put(18,32){\line(3,1){9.00}}
\qbezier(36,32)(34.5,30.5)(33,29)
\put(36,32){\line(3,1){9.00}}
\put(36,32){\line(-3,1){9.00}}
\qbezier(54,32)(52.5,30.5)(51,29)
\put(54,32){\line(3,1){9.00}}
\put(54,32){\line(-3,1){9.00}}
\qbezier(72,32)(70.5,30.5)(69,29)
\put(72,32){\line(-3,1){9.00}}


\put(12,40){\circle*{1.5}}
\put(30,40){\circle*{1.5}}
\put(48,40){\circle*{1.5}}
\put(3,43){\circle{1.5}}
\put(21,43){\circle{1.5}}
\put(39,43){\circle{1.5}}
\put(57,43){\circle{1.5}}
\put(6,46){\circle*{1.5}}
\put(24,46){\circle*{1.5}}
\put(42,46){\circle*{1.5}}
\put(60,46){\circle*{1.5}}
\put(15,49){\circle{1.5}}
\put(33,49){\circle{1.5}}
\put(51,49){\circle{1.5}}
\put(69,49){\circle{1.5}}
\put(18,52){\circle*{1.5}}
\put(36,52){\circle*{1.5}}
\put(54,52){\circle*{1.5}}
\put(72,52){\circle*{1.5}}
\put(27,55){\circle{1.5}}
\put(45,55){\circle{1.5}}
\put(63,55){\circle{1.5}}

\put(12,40){\line(-3,1){9.00}}
\put(12,40){\line(3,1){9.00}}
\put(30,40){\line(-3,1){9.00}}
\put(30,40){\line(3,1){9.00}}
\put(48,40){\line(-3,1){9.00}}
\put(48,40){\line(3,1){9.00}}

\qbezier(6,46)(4.5,44.5)(3,43)
\put(6,46){\line(3,1){9.00}}
\qbezier(24,46)(22.5,44.5)(21,43)
\put(24,46){\line(3,1){9.00}}
\put(24,46){\line(-3,1){9.00}}
\qbezier(42,46)(40.5,44.5)(39,43)
\put(42,46){\line(3,1){9.00}}
\put(42,46){\line(-3,1){9.00}}
\qbezier(60,46)(58.5,44.5)(57,43)
\put(60,46){\line(3,1){9.00}}
\put(60,46){\line(-3,1){9.00}}

\qbezier(18,52)(16.5,50.5)(15,49)
\put(18,52){\line(3,1){9.00}}
\qbezier(36,52)(34.5,50.5)(33,49)
\put(36,52){\line(3,1){9.00}}
\put(36,52){\line(-3,1){9.00}}
\qbezier(54,52)(52.5,50.5)(51,49)
\put(54,52){\line(3,1){9.00}}
\put(54,52){\line(-3,1){9.00}}
\qbezier(72,52)(70.5,50.5)(69,49)
\put(72,52){\line(-3,1){9.00}}
\color{red}
\qbezier(21,3)(16.5,11.5)(12,20)
\qbezier(39,3)(34.5,11.5)(30,20)
\qbezier(57,3)(52.5,11.5)(48,20)
\qbezier(15,9)(10.5,17.5)(6,26)
\qbezier(33,9)(28.5,17.5)(24,26)
\qbezier(51,9)(46.5,17.5)(42,26)
\qbezier(69,9)(64.5,17.5)(60,26)
\qbezier(27,15)(22.5,23.5)(18,32)
\qbezier(45,15)(40.5,23.5)(36,32)
\qbezier(63,15)(58.5,23.5)(54,32)

\qbezier(21,23)(16.5,31.5)(12,40)
\qbezier(39,23)(34.5,31.5)(30,40)
\qbezier(57,23)(52.5,31.5)(48,40)
\qbezier(15,29)(10.5,37.5)(6,46)
\qbezier(33,29)(28.5,37.5)(24,46)
\qbezier(51,29)(46.5,37.5)(42,46)
\qbezier(69,29)(64.5,37.5)(60,46)
\qbezier(27,35)(22.5,43.5)(18,52)
\qbezier(45,35)(40.5,43.5)(36,52)
\qbezier(63,35)(58.5,43.5)(54,52)
\qbezier(12,0)(12.5,-1)(13,-2)
\qbezier(30,0)(30.5,-1)(31,-2)
\qbezier(48,0)(48.5,-1)(49,-2)
\qbezier(6,6)(6.5,5)(7,4)
\qbezier(24,6)(24.5,5)(25,4)
\qbezier(42,6)(42.5,5)(43,4)
\qbezier(60,6)(60.5,5)(61,4)
\qbezier(18,12)(18.5,11)(19,10)
\qbezier(36,12)(36.5,11)(37,10)
\qbezier(54,12)(54.5,11)(55,10)
\qbezier(72,12)(72.5,11)(73,10)
\qbezier(72,32)(72.5,31)(73,30)
\qbezier(72,52)(72.5,51)(73,50)

\qbezier(3,3)(2.5,4)(2,5)
\qbezier(3,23)(2.5,24)(2,25)
\qbezier(3,43)(2.5,44)(2,45)
\qbezier(21,43)(20.5,44)(20,45)
\qbezier(39,43)(38.5,44)(38,45)
\qbezier(57,43)(56.5,44)(56,45)
\qbezier(15,49)(14.5,50)(14,51)
\qbezier(33,49)(32.5,50)(32,51)
\qbezier(51,49)(50.5,50)(50,51)
\qbezier(69,49)(68.5,50)(68,51)

\qbezier(27,55)(26.5,56)(26,57)
\qbezier(45,55)(44.5,56)(44,57)
\qbezier(63,55)(62.5,56)(62,57)
\put(15.3,15){$\be_o$}
\put(15.2,14){\vector(-1,2){1}}
\end{picture}
\vspace{5mm}
\caption{\scriptsize \emph{Step 3.} \emph{a}) Copies of the unperturbed graph $\cG$ stacked together along the vector $\ga_3$.  \emph{b}) The limit graph $\wt\cG$ is generated by infinitely many copies of $\cG$. The copies are connected by the red edges $\be_o+\ga$, $\forall\,\ga\in\wt\G$, where $\be_o=(v_1,v_2+\ga_3)$ and $\wt\G$ is the lattice with basis $\ga_1,\ga_2,\ga_3$.} \label{figS3}
\end{figure}
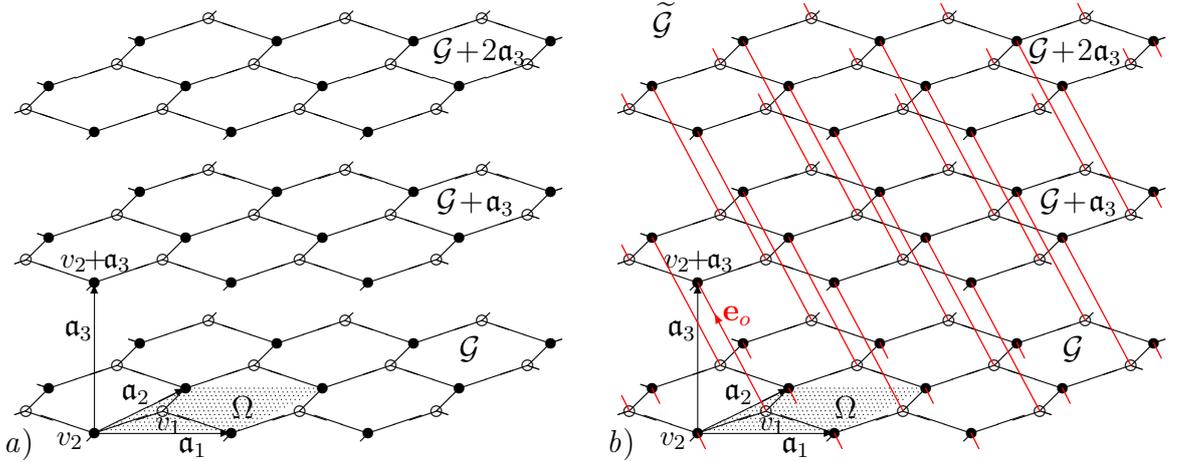

\begin{remark}
\emph{i}) The perturbed graph $\cG_t$ and its limit graph $\wt\cG$ have the same fundamental graph $\wt\cG_*=(\cV_*,\wt\cE_*)$ given by (\ref{pefg}) but with distinct edge indices denoted by $\t:\wt\cA_*\to\Z^{d}$ and $\wt\t:\wt\cA_*\to\Z^{d+1}$, respectively. Here $\wt\cA_*$ is the set of doubled oriented edges of $\wt\cG_*$. The relation between $\t$ and $\wt\t$ is given by \er{wtta}.

\emph{ii}) The limit graph $\wt\cG$ constructed above is actually an infinite-fold cover of the perturbed graph $\cG_t$ (see, e.g., \cite{S13} for more details on covering graphs).
\end{remark}

For readers convenience, in the following table we summarize different types of graphs we deal with.

\begin{table}[h]
\begin{tabular}{|l|l|c|}
  \hline
\rule{0pt}{12pt}  &\hspace{20mm}\emph{Periodic graph} & \emph{Fundamental graph} \\ \hline
\textbf{Unperturbed} \rule{0pt}{12pt} &$\cG=(\cV,\cE)\ss\R^d$ is $\G$-periodic & $\cG_*=\cG/\G=(\cV_*,\cE_*)$  \\
\textbf{graph $\cG$}& &  \\\hline
\textbf{Perturbed} \rule{0pt}{12pt}  & $\cG_t=(\cV,\cE_t)\ss\R^d$ is $\G$-periodic, & \\
\textbf{graph $\cG_t$}    & $\cE_t=\cE\cup\{\be_o+\ga : \ga\in\G\}$, & \\
$t=(t_1,\ldots,t_d)\in\Z^d$  & where  $\be_o=(v_1,v_2+t_1\ga_1+\ldots+t_d\ga_d)$ & $\wt\cG_*=(\cV_*,\wt\cE_*)$ \\
    & for some fixed $v_1,v_2\in\cV\cap\Omega$ & $\wt\cE_*=\cE_*\cup\{\be_o\}$\\\cline{1-2}
\textbf{Limit} \rule{0pt}{14pt} & $\wt\cG=(\wt\cV,\wt\cE\,)\ss\R^{d+1}$ is $\wt\G$-periodic,  & $\be_o=(v_1,v_2)$ \\
\textbf{graph $\wt\cG$}& $\wt\cV=\bigcup\limits_{n\in\Z}(\cV+n\ga_{d+1})$,  & $v_1,v_2\in\cV_*$\\ & $\wt\cE=\Big(\bigcup\limits_{n\in\Z}(\cE+n\ga_{d+1})\Big)\cup\{\be_o+\ga :\ga\in\wt\G\}$, & \\
& where $\be_o=(v_1,v_2+\ga_{d+1})$  & \\
\hline
\multicolumn{3}{|l|}{$\G\ss\R^d$ \rule{0pt}{12pt}  is a lattice with basis $\ga_1,\ldots,\ga_d$ and the fundamental cell $\Omega$}\\\hline
\multicolumn{3}{|l|}{$\wt\G\ss\R^{d+1}$ \rule{0pt}{14pt}  is a lattice with basis $\ga_1,\ldots,\ga_d,\ga_{d+1}$}\\\hline
\end{tabular}
\end{table}

\subsection{Asymptotically isospectral graphs}
Let $H_{\cG_t}=\D_{\cG_t}+Q$ be the Schr\"odinger operator with a periodic potential $Q$ on the perturbed graph $\cG_t=(\cV,\cE_t)$ with the fundamental graph $\wt\cG_*=(\cV_*,\wt\cE_*)$. The periodic potential $Q:\cV\to\R$ is uniquely defined by the vector
\[\lb{depo}
(Q_v)_{v\in\cV_*}\in\R^\n,\qqq\textrm{where}\qqq \n=\#\cV_*.
\]
Recall that the limit graph $\wt\cG$ has the same fundamental graph $\wt\cG_*=(\cV_*,\wt\cE_*)$.

We consider the Schr\"odinger operator $H_{\wt\cG}=\D_{\wt\cG}+Q$ on the limit graph $\wt\cG=(\wt\cV,\wt\cE\,)$ with the periodic potential $Q:\wt\cV\to\R$ defined by the vector \er{depo}. Due to (\ref{ban.1H}), (\ref{specH}), the spectrum of $H_{\wt\cG}$ has the form
\[\lb{debe}
\begin{array}{l}
\displaystyle \s(H_{\wt\cG})=
\bigcup_{j=1}^\n\s_j(H_{\wt\cG}),\qqq
\s_j(H_{\wt\cG})=
\big[\,\wt\l_j^-,\wt\l_j^+\big],\qqq \n=\#\cV_*,\\
\displaystyle\textrm{where}\qqq \wt\l_j^-=\min_{k\in\T^{d+1}}\wt\l_j(k),\qqq \wt\l_j^+=\max_{k\in\T^{d+1}}\wt\l_j(k),
\end{array}
\]
and $\wt\l_j(\cdot)$ are the band functions of $H_{\wt\cG}$.

In the following theorem we show that the Schr\"odinger operators $H_{\cG_t}$ and $H_{\wt\cG}$ have asymptotically close spectra.

\begin{theorem}\lb{TAig}
The perturbed graph $\cG_t$ is asymptotically isospectral to its limit graph $\wt\cG$, i.e.,
\[\lb{as00}
\wt\l_j^-=\lim\limits_{|t|\to\iy}\l_j^-(t),\qquad
\wt\l_j^+=\lim\limits_{|t|\to\iy}\l_j^+(t),\qqq \forall\,j\in\N_\n,\qqq \n=\#\cV_*,
\]
where $\wt\l_j^\pm$ and $\l_j^\pm(t)$ are the band edges of the Schr\"odinger operators $H_{\wt\cG}=\D_{\wt\cG}+Q$ on $\wt\cG$ and $H_{\cG_t}=\D_{\cG_t}+Q$ on $\cG_t$, respectively.
\end{theorem}

\begin{remark}
\emph{i}) The perturbed graph $\cG_t$ has dimension $d$, but its limit graph $\wt\cG$ has dimension $d+1$. Thus, adding sufficiently long edges (in a periodic way) to a periodic graph leads to an asymptotically isospectral periodic graph of a higher dimension but without long edges. This effect might have the following explanation. A perturbation by a long edge means that its endpoints which are distant from each other on the unperturbed graph should interact on the perturbed graph, i.e., should become close to each other. But the only possibility to make distant vertices close is to embed the graph in the space of a higher dimension.

\emph{ii}) Due to Theorem \ref{TAig}, the spectrum of the Schr\"odinger operator on the limit graph $\wt\cG$ can be approximated by spectra of the Schr\"odinger operators on its lower dimensional subcovering graphs $\cG_t$.
\end{remark}

Now we formulate a criterion for the perturbed graph $\cG_t$ to be not only asymptotically isospectral but just isospectral to its limit graph $\wt\cG$. Denote by $K_j^\pm$ the level sets corresponding to the band edges $\wt\l_j^\pm$ of the Schr\"odinger operator $H_{\wt\cG}$ on the limit graph $\wt\cG$, i.e.,
\[\lb{Kjpm}
K_j^\pm=\big\{k\in\T^{d+1} : \wt\l_j(k)=\wt\l_j^\pm\big\},\qqq j\in\N_\n, \qqq \n=\#\cV_*,
\]
where $\wt\l_j(\cdot)$ are the band functions of $H_{\wt\cG}$.

\begin{theorem}\lb{TJig}
Let $t=(t_1,\ldots,t_d)\in\Z^d$ and $j\in\N_\n$. Then the band edge $\l_j^\pm(t)$ of the Schr\"odinger operator $H_{\cG_t}$ on the perturbed graph $\cG_t$ coincides with the corresponding band edge $\wt\l_j^\pm$ of $H_{\wt\cG}$ on its limit graph $\wt\cG$, i.e., $\l_j^\pm(t)=\wt\l_j^\pm$, if and only if
\[\lb{soeq}
t_1k_1^o+\ldots+t_dk_d^o  \equiv k_{d+1}^o\; (\hspace{-4mm}\mod 2\pi)
\]
for some $k_o=(k_1^o,\ldots,k_d^o,k_{d+1}^o)\in K_j^\pm$, where $K_j^\pm$ is given by \er{Kjpm}. In particular, if $0\in K_j^\pm$, then $\l_j^\pm(t)=\wt\l_j^\pm$.
\end{theorem}

The following criterion of isospectrality is a direct consequence of Theorem \ref{TJig}.

\begin{corollary}\lb{CJig}
Let $t=(t_1,\ldots,t_d)\in\Z^d$. Then the perturbed graph $\cG_t$ is isospectral to its limit graph $\wt\cG$, i.e., $\s(H_{\cG_t})=\s(H_{\wt\cG})$ if and only if for each $j\in\N_\n$ the condition \er{soeq} holds true for some $k^-\in K_j^-$ and for some $k^+\in K_j^+$, where $K_j^\pm$ are given by \er{Kjpm}.
\end{corollary}

\subsection{Asymptotics of the band edges.} In this subsection we present asymptotics of the band edges $\l_j^\pm(t)$ of the Schr\"odinger operator $H_{\cG_t}$ on the perturbed graph $\cG_t$ as $|t|\to\iy$.

For each $t=(t_1,\ldots,t_d)\in\Z^d$ we define the function $\vr_t:\T^{d+1}\to\R$ by
\[\lb{fual}
\textstyle
\vr_t(k)=\frac1{2\pi}(t_1k_1+\ldots+t_dk_d-k_{d+1}),\qqq k=(k_1,\ldots,k_d,k_{d+1})\in\T^{d+1}.
\]
Denote by $\a(x)$ the difference between $x\in\R$ and the integer $\mn(x)$ closest to $x$, i.e.,
\[\lb{alox}
\textstyle\a(x)=x-\mn(x),\qqq \a(x)\in\big(-\frac12,\frac12\big],\qqq \forall\, x\in\R.
\]

From now on we impose the following assumptions on a band function $\wt\l_j(\cdot)$ of the Schr\"odinger operator $H_{\wt\cG}$ on the limit graph $\wt\cG$:

\medskip

\textbf{Assumption A}

\textbf{A1} The band function $\wt\l_j(\cdot)$ of $H_{\wt\cG}$ has a non-degenerate minimum $\wt\l_j^-$ (maximum $\wt\l_j^+$) at some point $k_o\in\T^{d+1}$, and this extremum $\wt\l_j^\pm$ is attained by the single band function $\wt\l_j(\cdot)$.

\textbf{A2} $k_o$ is the only (up to evenness) minimum (maximum) point of $\wt\l_j(\cdot)$ in $\T^{d+1}$, see Remark~\ref{prbf}.\emph{ii}.

\begin{theorem}\lb{TAsy}
Let $t=(t_1,\ldots,t_d)\in\Z^d$ and let for some $j\in\N_\n$ the band function $\wt\l_j(\cdot)$ of the Schr\"odinger operator $H_{\wt\cG}$ on the limit graph $\wt\cG$ satisfy the Assumption A. Then the lower band edge $\l_j^-(t)$ (the upper band edge $\l_j^+(t)$) of the Schr\"odinger operator $H_{\cG_t}$ on the perturbed graph $\cG_t$ has the following asymptotics
\[\lb{aslk}
\begin{array}{l}
\textstyle\l_j^\pm(t)=\wt\l_j^\pm\mp\dfrac{2\pi^2\a^2(\vr_t(k_o))}{|\bH^{-1/2}\gt|^2}+g(t),\qq \textrm{where}\qq  g(t)=O(\frac1{|t|^3}),\qq \textrm{as}\qq |t|\to\iy,\\[12pt]
\bH=\mp\mathrm{Hess}\,\wt\l_j(k_o),\qqq \gt=(t_1,\ldots,t_d,-1),
\end{array}
\]
$\a(\cdot)$ and $\vr_t(\cdot)$ are given by \er{alox} and \er{fual}, respectively; $\mathrm{Hess}\,\wt\l_j(k_o)$ is the Hessian of $\wt\l_j$ at $k_o$. Moreover, if all components of $k_o$ are equal to either 0 or $\pi$, then $g(t)=O(\frac1{|t|^4})$.
\end{theorem}

\begin{remark}\lb{re26}
\emph{i}) Assumption A1 is often used to establish many important properties of periodic operators such as electron's effective masses \cite{AM76}, Green's function  asymptotics \cite{KKR17,K18} and other issues. Although in the continuous case it is commonly believed that the band edges are non-degenerate for "generic" potentials, in discrete settings there are known counterexamples \cite{FK18,BK20}, as well as various positive results in some cases \cite{DKS20}.

\emph{ii}) The asymptotics \er{aslk} agrees with Theorem~\ref{TJig}. Indeed,
if for $t\in\Z^d$ the identity \er{soeq} holds, then, due to Theorem~\ref{TJig}, $\l_j^\pm(t)=\wt\l_j^\pm$. On the other hand, for such $t$ we have $\vr_t(k_o)\in\Z$, $\a(\vr_t(k_o))=0$, and the asymptotics \er{aslk} has the form $\l_j^\pm(t)=\wt\l_j^\pm+O(\frac1{|t|^3})$.

\emph{iii}) The asymptotics \er{aslk} can be carried over to the case when the band edge $\wt\l_j^\pm$ occurs at finitely many quasimomenta $k_o\in\T^{d+1}$ (instead of assuming the condition A2) by taking minimum (maximum) among the asymptotics coming from all these non-degenerate isolated extrema.

\emph{iv}) For each $t\in\Z^d$ the spectral bands
$$
\s_{j}(H_{\cG_t})=\big[\l_j^-(t),\l_j^+(t)\big]\qqq\textrm{and}\qqq \s_j(H_{\wt\cG})=\big[\,\wt\l_j^-,\wt\l_j^+\big],\qqq j\in\N_\n, \qqq \n=\#\cV_*,
$$
of the Schr\"odinger operators $H_{\cG_t}$ and $H_{\wt\cG}$, respectively, satisfy $\s_{j}(H_{\cG_t})\subseteq\s_j(H_{\wt\cG})$ (see Proposition \ref{Pnin0}). In particular, if $\s_j(H_{\wt\cG})$ is flat, then $\s_{j}(H_{\cG_t})$ is also flat for any $t\in\Z^d$.
\end{remark}

\subsection{Examples.} In this subsection we present some simple examples of asymptotically isospectral and just isospectral periodic graphs. The proofs of the examples are based on Theorems \ref{TAig}, \ref{TJig}, \ref{TAsy} and given in Section \ref{Sec4}.

\medskip

First we consider the Schr\"odinger operator $H_\dL=\D_\dL+Q$ with a 2-periodic potential $Q$ on the one-dimensional lattice $\dL$, Fig.~\ref{fig5}\emph{a}. Without loss of generality we may assume that
\[\lb{poe2}
Q_0=-Q_1=q>0.
\]
Then the spectrum of $H_\dL$ has the form
\[\lb{spHL}
\s(H_\dL)=\big[2-\sqrt{4+q^2}\,,2-q\big]\cup\big[2+q,2+\sqrt{4+q^2}\,\big]
\]
(see Section \ref{Ss4.2}). We describe the spectrum of the Schr\"odinger operator $H_{\dL_t}=\D_{\dL_t}+Q$ with the same potential $Q$ on the graph $\dL_t$, $t\in\Z$, shown in Fig.~\ref{fig5}\emph{b}.

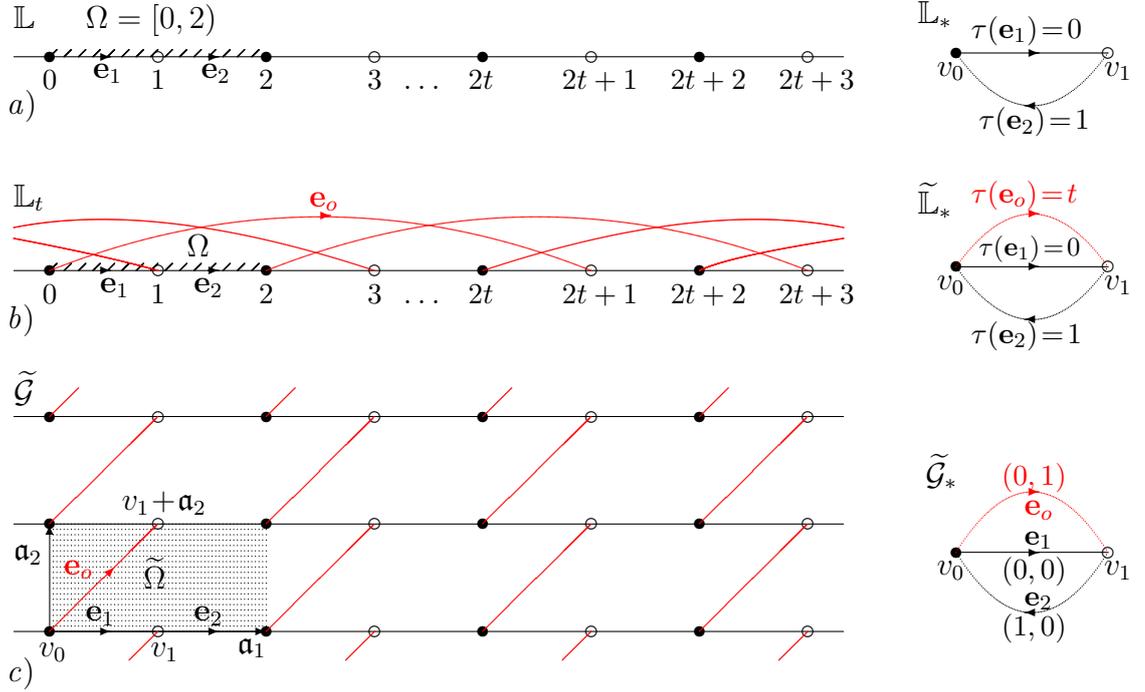
\begin{figure}[t!]\centering
\unitlength 0.95mm 
\begin{picture}(125,13)
\hspace{-5mm}
\put(5,10){\line(1,0){115.00}}
\bezier{60}(10,10)(10.5,10.5)(11,11)
\bezier{60}(12,10)(12.5,10.5)(13,11)
\bezier{60}(14,10)(14.5,10.5)(15,11)
\bezier{60}(16,10)(16.5,10.5)(17,11)
\bezier{60}(18,10)(18.5,10.5)(19,11)
\bezier{60}(20,10)(20.5,10.5)(21,11)
\bezier{60}(22,10)(22.5,10.5)(23,11)
\bezier{60}(24,10)(24.5,10.5)(25,11)
\bezier{60}(26,10)(26.5,10.5)(27,11)
\bezier{60}(28,10)(28.5,10.5)(29,11)
\bezier{60}(30,10)(30.5,10.5)(31,11)
\bezier{60}(32,10)(32.5,10.5)(33,11)
\bezier{60}(34,10)(34.5,10.5)(35,11)
\bezier{60}(36,10)(36.5,10.5)(37,11)
\bezier{60}(38,10)(38.5,10.5)(39,11)
\put(15,14){$\Omega=[0,2)$}
\put(5,14){$\dL$}
\put(32.5,10){\vector(1,0){1}}
\put(17.5,10){\vector(1,0){1}}
\put(16,7.0){$\be_1$}
\put(31,7.0){$\be_2$}

\put(10,10){\circle*{1.5}}
\put(25,10){\circle{1.5}} \put(40,10){\circle*{1.5}}
\put(55,10){\circle{1.5}} \put(70,10){\circle*{1.5}}
\put(85,10){\circle{1.5}} \put(100,10){\circle*{1.5}} \put(115,10){\circle{1.5}}
\put(9,5.5){\small$0$}
\put(24,5.5){\small$1$} \put(39,5.5){\small$2$} \put(54,5.5){\small$3$}\put(59,5.5){$\ldots$}
\put(68,5.5){\small$2t$} \put(81,5.5){\small$2t+1$} \put(96,5.5){\small$2t+2$}\put(111,5.5){\small$2t+3$}
\put(4,2){\emph{a})}
\end{picture}\unitlength 1.0mm
\begin{picture}(25,20)
\put(5,10){\line(1,0){20.00}}
\put(15,10){\vector(1,0){1.00}}
\put(15,3){\vector(-1,0){1.00}}
\put(5,10){\circle*{1.5}}
\put(25,10){\circle{1.5}}
\bezier{100}(5,10)(15,-4)(25,10)
\put(0,14){$\dL_*$}
\put(2.5,7){$v_0$}
\put(24.5,7){$v_1$}

\put(7,12){\small$\t(\be_1)\!=\!0$}
\put(8,-0.1){\small$\t(\be_2)\!=\!1$}
\end{picture}

\unitlength 0.95mm
\begin{picture}(125,25)
\hspace{-5mm}
\put(5,10){\line(1,0){115.00}} \put(10,10){\circle*{1.5}}
\put(25,10){\circle{1.5}} \put(40,10){\circle*{1.5}}
\put(55,10){\circle{1.5}} \put(70,10){\circle*{1.5}}
\put(85,10){\circle{1.5}} \put(100,10){\circle*{1.5}} \put(115,10){\circle{1.5}}
\put(9,5.5){\small$0$}
\put(24,5.5){\small$1$} \put(39,5.5){\small$2$} \put(54,5.5){\small$3$}\put(59,5.5){$\ldots$}
\put(68,5.5){\small$2t$}\put(81,5.5){\small$2t+1$} \put(96,5.5){\small$2t+2$}\put(111,5.5){\small$2t+3$}
\put(5,19){$\dL_t$}
\put(32.5,10){\vector(1,0){1}}
\put(17.5,10){\vector(1,0){1}}
\put(17,7.0){$\be_1$}
\put(30,7.0){$\be_2$}

\bezier{60}(10,10)(10.5,10.5)(11,11)
\bezier{60}(12,10)(12.5,10.5)(13,11)
\bezier{60}(14,10)(14.5,10.5)(15,11)
\bezier{60}(16,10)(16.5,10.5)(17,11)
\bezier{60}(18,10)(18.5,10.5)(19,11)
\bezier{60}(20,10)(20.5,10.5)(21,11)
\bezier{60}(22,10)(22.5,10.5)(23,11)
\bezier{60}(24,10)(24.5,10.5)(25,11)
\bezier{60}(26,10)(26.5,10.5)(27,11)
\bezier{60}(28,10)(28.5,10.5)(29,11)
\bezier{60}(30,10)(30.5,10.5)(31,11)
\bezier{60}(32,10)(32.5,10.5)(33,11)
\bezier{60}(34,10)(34.5,10.5)(35,11)
\bezier{60}(36,10)(36.5,10.5)(37,11)
\bezier{60}(38,10)(38.5,10.5)(39,11)

\put(29,12){$\Omega$}

\put(4,2){\emph{b})}
\color{red}
\bezier{600}(10,10)(47.5,25)(85,10)
\bezier{600}(40,10)(77.5,25)(115,10)
\bezier{600}(70,10)(100,20)(120,16)
\bezier{600}(100,10)(104,12)(120,14.5)
\bezier{600}(25,10)(18,12)(5,14.5)
\bezier{600}(55,10)(25,20)(5,16)

\put(46,19){$\be_o$}
\put(48,17.7){\vector(1,0){1}}
\end{picture}\unitlength 1.0mm
\begin{picture}(25,28)
\put(5,10){\line(1,0){20.00}}
\put(15,3){\vector(-1,0){1.00}}
\put(15,10){\vector(1,0){1.00}}
\put(5,10){\circle*{1.5}}
\put(25,10){\circle{1.5}}
\bezier{100}(5,10)(15,-4)(25,10)
\put(0,17){$\wt\dL_*$}
\put(2.5,7){$v_0$}
\put(24.5,7){$v_1$}
\put(8,11.5){\footnotesize$\t(\be_1)\!=\!0$}
\put(7,-0.1){\small$\t(\be_2)\!=\!1$}
\color{red}
\bezier{100}(5,10)(15,24)(25,10)
\put(15,17){\vector(1,0){1.00}}
\put(7,18.5){\small$\t(\be_o)\!=\!t$}
\end{picture}

\unitlength 0.95mm
\begin{picture}(125,50)
\hspace{-5mm}
\put(5,10){\line(1,0){115.00}}
\put(5,25){\line(1,0){115.00}}
\put(5,40){\line(1,0){115.00}}
\put(32.5,10){\vector(1,0){1}}
\put(17.5,10){\vector(1,0){1}}

\put(15,11.5){$\be_1$}
\put(30,11.5){$\be_2$}

\put(23,16){$\wt\Omega$}
\put(5,42){$\wt\cG$}
\put(10,10){\vector(1,0){30}}
\put(10,10){\vector(0,1){15}}
\put(36,7.0){$\ga_1$}
\put(5,20){$\ga_2$}

\put(8.5,6.5){$v_0$}
\put(24,6.5){$v_1$}
\put(20,27.0){$v_1\!+\!\ga_2$}

\bezier{20}(10.5,10)(10.5,17.5)(10.5,25)
\bezier{20}(11.0,10)(11.0,17.5)(11.0,25)
\bezier{20}(11.5,10)(11.5,17.5)(11.5,25)
\bezier{20}(12.0,10)(12.0,17.5)(12.0,25)
\bezier{20}(12.5,10)(12.5,17.5)(12.5,25)
\bezier{20}(13.0,10)(13.0,17.5)(13.0,25)
\bezier{20}(13.5,10)(13.5,17.5)(13.5,25)
\bezier{20}(14.0,10)(14.0,17.5)(14.0,25)
\bezier{20}(14.5,10)(14.5,17.5)(14.5,25)
\bezier{20}(15.0,10)(15.0,17.5)(15.0,25)
\bezier{20}(15.5,10)(15.5,17.5)(15.5,25)
\bezier{20}(16.0,10)(16.0,17.5)(16.0,25)
\bezier{20}(16.5,10)(16.5,17.5)(16.5,25)
\bezier{20}(17.0,10)(17.0,17.5)(17.0,25)
\bezier{20}(17.5,10)(17.5,17.5)(17.5,25)
\bezier{20}(18.0,10)(18.0,17.5)(18.0,25)
\bezier{20}(18.5,10)(18.5,17.5)(18.5,25)
\bezier{20}(19.0,10)(19.0,17.5)(19.0,25)
\bezier{20}(19.5,10)(19.5,17.5)(19.5,25)
\bezier{20}(20.0,10)(20.0,17.5)(20.0,25)
\bezier{20}(20.5,10)(20.5,17.5)(20.5,25)
\bezier{20}(21.0,10)(21.0,17.5)(21.0,25)
\bezier{20}(21.5,10)(21.5,17.5)(21.5,25)
\bezier{20}(22.0,10)(22.0,17.5)(22.0,25)
\bezier{20}(22.5,10)(22.5,17.5)(22.5,25)
\bezier{20}(23.0,10)(23.0,17.5)(23.0,25)
\bezier{20}(23.5,10)(23.5,17.5)(23.5,25)
\bezier{20}(24.0,10)(24.0,17.5)(24.0,25)
\bezier{20}(24.5,10)(24.5,17.5)(24.5,25)
\bezier{20}(25.0,10)(25.0,17.5)(25.0,25)
\bezier{20}(25.5,10)(25.5,17.5)(25.5,25)
\bezier{20}(26.0,10)(26.0,17.5)(26.0,25)
\bezier{20}(26.5,10)(26.5,17.5)(26.5,25)
\bezier{20}(27.0,10)(27.0,17.5)(27.0,25)
\bezier{20}(27.5,10)(27.5,17.5)(27.5,25)
\bezier{20}(28.0,10)(28.0,17.5)(28.0,25)
\bezier{20}(28.5,10)(28.5,17.5)(28.5,25)
\bezier{20}(29.0,10)(29.0,17.5)(29.0,25)
\bezier{20}(29.5,10)(29.5,17.5)(29.5,25)
\bezier{20}(30.0,10)(30.0,17.5)(30.0,25)
\bezier{20}(30.5,10)(30.5,17.5)(30.5,25)
\bezier{20}(31.0,10)(31.0,17.5)(31.0,25)
\bezier{20}(31.5,10)(31.5,17.5)(31.5,25)
\bezier{20}(32.0,10)(32.0,17.5)(32.0,25)
\bezier{20}(32.5,10)(32.5,17.5)(32.5,25)
\bezier{20}(33.0,10)(33.0,17.5)(33.0,25)
\bezier{20}(33.5,10)(33.5,17.5)(33.5,25)
\bezier{20}(34.0,10)(34.0,17.5)(34.0,25)
\bezier{20}(34.5,10)(34.5,17.5)(34.5,25)
\bezier{20}(35.0,10)(35.0,17.5)(35.0,25)
\bezier{20}(35.5,10)(35.5,17.5)(35.5,25)
\bezier{20}(36.0,10)(36.0,17.5)(36.0,25)
\bezier{20}(36.5,10)(36.5,17.5)(36.5,25)
\bezier{20}(37.0,10)(37.0,17.5)(37.0,25)
\bezier{20}(37.5,10)(37.5,17.5)(37.5,25)
\bezier{20}(38.0,10)(38.0,17.5)(38.0,25)
\bezier{20}(38.5,10)(38.5,17.5)(38.5,25)
\bezier{20}(39.0,10)(39.0,17.5)(39.0,25)
\bezier{20}(39.5,10)(39.5,17.5)(39.5,25)
\bezier{20}(40.0,10)(40.0,17.5)(40.0,25)

\put(10,10){\circle*{1.5}}
\put(25,10){\circle{1.5}} \put(40,10){\circle*{1.5}}
\put(55,10){\circle{1.5}} \put(70,10){\circle*{1.5}}
\put(85,10){\circle{1.5}} \put(100,10){\circle*{1.5}} \put(115,10){\circle{1.5}}

\put(10,25){\circle*{1.5}}
\put(25,25){\circle{1.5}} \put(40,25){\circle*{1.5}}
\put(55,25){\circle{1.5}} \put(70,25){\circle*{1.5}}
\put(85,25){\circle{1.5}} \put(100,25){\circle*{1.5}} \put(115,25){\circle{1.5}}

\put(10,40){\circle*{1.5}}
\put(25,40){\circle{1.5}} \put(40,40){\circle*{1.5}}
\put(55,40){\circle{1.5}} \put(70,40){\circle*{1.5}}
\put(85,40){\circle{1.5}} \put(100,40){\circle*{1.5}} \put(115,40){\circle{1.5}}
\put(4,3){\emph{c})}
\color{red}

\put(10,10){\line(1,1){15}}
\put(40,10){\line(1,1){15}}
\put(70,10){\line(1,1){15}}
\put(100,10){\line(1,1){15}}

\put(10,25){\line(1,1){15}}
\put(40,25){\line(1,1){15}}
\put(70,25){\line(1,1){15}}
\put(100,25){\line(1,1){15}}

\bezier{60}(10,40)(12,42)(14,44)
\bezier{60}(40,40)(42,42)(44,44)
\bezier{60}(70,40)(72,42)(74,44)
\bezier{60}(100,40)(102,42)(104,44)

\bezier{60}(25,10)(23,8)(21,6)
\bezier{60}(55,10)(53,8)(51,6)
\bezier{60}(85,10)(83,8)(81,6)
\bezier{60}(115,10)(113,8)(111,6)

\put(12,18){$\be_o$}
\put(18,18){\vector(1,1){1.00}}
\end{picture}\unitlength 1.0mm
\begin{picture}(25,40)
\put(5,20){\line(1,0){20.00}}
\put(15,20){\vector(1,0){1.00}}
\put(15,12){\vector(-1,0){1.00}}
\put(5,20){\circle*{1.5}}
\put(25,20){\circle{1.5}}
\bezier{100}(5,20)(15,4)(25,20)
\put(1,29){$\wt\cG_*$}
\put(2.5,17){$v_0$}
\put(24.5,17){$v_1$}

\put(14,21.0){\small$\be_1$}
\put(11,16.5){\small$(0,0)$}
\put(11,9){\small$(1,0)$}
\put(14,13){\small$\be_2$}
\color{red}
\put(15,28){\vector(1,0){1.00}}
\bezier{100}(5,20)(15,36)(25,20)

\put(14,25){$\be_o$}
\put(11,29){\small$(0,1)$}
\end{picture}

\caption{\scriptsize\emph{a}) The one-dimensional lattice $\dL$. \emph{b}) The perturbed lattice $\dL_t$ with period 2; the fundamental cell $\Omega=[0,2)$; the added edges have red color. \emph{c}) The limit graph $\wt\cG$ of $\dL_t$ is obtained by stacking together infinitely many black copies of $\dL$ along the vector $\ga_2$. The copies are connected by the red edges $\be_o+\ga$, $\forall\,\ga\in\wt\G$, where $\be_o=(v_0,v_1+\ga_2)$, and $\wt\G$ is the lattice with basis $\ga_1,\ga_2$ and the fundamental cell $\wt\Omega$.  \emph{a}) -- \emph{c}) (right) The fundamental graphs $\dL_*=\dL/2\Z$ of $\dL$; $\wt\dL_*$ of $\dL_t$ and $\wt\cG_*$ of $\wt\cG$; the indices are shown near the edges.} \label{fig5}
\end{figure}

\begin{example}\lb{Ehex} Let $t\in\mathbb{Z}$, and let $\dL_t$ be the perturbed graph obtained from the lattice $\dL$ by adding an edge $\mathbf{e}_o$ with index $\tau(\mathbf{e}_o)=t$ to its fundamental graph $\dL_*=\dL/2\Z$, see Fig.~\ref{fig5}\emph{b}. Then the perturbed lattice $\dL_t$ is asymptotically isospectral to the graph $\wt\cG$ shown in Fig.~\ref{fig5}\emph{c} as $|t|\to\iy$, and the spectrum $\s(H_{\dL_t})$ of the Schr\"odinger operator $H_{\dL_t}=\D_{\dL_t}+Q$ with the 2-periodic potential $Q$ defined by \er{poe2} on $\dL_t$ has the form
\[\lb{stin}
\s(H_{\dL_t})=\big[3-\sqrt{9+q^2}\,,3-q-\ve(t)\big]\cup
\big[3+q+\ve(t),3+\sqrt{9+q^2}\,\big].
\]

\emph{\,i}) If $t-1\in3\Z$, then $\ve(t)=0$, and $\s(H_{\dL_t})=\s(H_{\wt\cG})$.

\emph{ii}) If $t-1\notin3\Z$, then
$$
\textstyle\ve(t)=\frac{\pi^2}{6q}\cdot\frac1{t^2}+O(\frac1{|t|^3}),\qqq \textrm{as} \qqq|t|\to\infty.
$$
\end{example}

\begin{remark} \emph{i}) The graph $\wt\cG$ shown in Fig.~\ref{fig5}\emph{c} is isomorphic to the hexagonal lattice (see Fig.~\ref{figS1}\emph{a}). Thus, the perturbed lattice $\dL_t$ is asymptotically isospectral to the hexagonal lattice as $|t|\to\iy$.

\emph{ii}) If $t-1\in3\Z$, then the Schr\"odinger operators on the perturbed lattice $\dL_t$ and on the graph $\wt\cG$ shown in Fig.~\ref{fig5}\emph{c} have the same spectrum, i.e., the perturbed lattice $\dL_t$ and the graph $\wt\cG$ are isospectral.

\emph{iii}) For any $t\in\Z$, the spectrum of the Schr\"odinger operator $H_{\dL_t}$ on the perturbed lattice $\dL_t$ is symmetric with respect to the point 3 and has a gap of length $2(q+\ve(t))>0$.
\end{remark}

As the second example we consider the $d$-dimensional lattice $\dL^d$, for $d=2,3$ see Fig.~\ref{fig2}\emph{a,c}. The spectrum of the Laplacian $\D_{\dL^d}$ on $\dL^d$ has the form $\s(\D_{\dL^d})=[0,4d]$. We describe the Laplacian spectrum under the perturbations of the lattice by adding an edge to its fundamental graph.

\begin{example}\lb{ExSl} Let $t=(t_1,\ldots,t_d)\in\Z^d$, and let $\dL^d_t$ be the perturbed graph obtained from the lattice $\dL^d$ by adding an edge $\be_o$ with index $\t(\be_o)=t$ to its fundamental graph $\dL^d_*$ (for $d=2$ see Fig.~\ref{fig2}\emph{b}). Then the perturbed $d$-dimensional lattice $\dL^d_t$ is asymptotically isospectral to the $(d+1)$-dimensional lattice $\dL^{d+1}$ as $|t|\to\iy$ and the spectrum $\s(\D_{\dL^d_t})$ of the Laplacian $\D_{\dL^d_t}$ on $\dL^d_t$ is given by
\[\lb{sDLt}
\s(\D_{\dL^d_t})=[0,\l_1^+(t)].
\]

\emph{\,i}) If $t_1+\ldots+t_d$ is odd, then $\l_1^+(t)=4(d+1)$ and $\s(\D_{\dL^d_t})=\s(\D_{\dL^{d+1}})$.

\emph{ii}) If $t_1+\ldots+t_d$ is even, then
\[\lb{ase1}
\textstyle\l_1^+(t)=4(d+1)-\frac{\pi^2}{|t|^2}+O\big(\frac1{|t|^4}\big),
\qqq \textrm{as}\qqq |t|\to\iy.
\]
\end{example}

\begin{remark} If $t_1+\ldots+t_d$ is odd, then $\l_1^+(t)$ does not depend on $t$ and the Laplacians on the perturbed $d$-dimensional lattice $\dL^d_t$ and on the $(d+1)$-dimensional lattice $\dL^{d+1}$ have the same spectrum, i.e., the lattices $\dL^d_t$ and $\dL^{d+1}$ are just isospectral.
\end{remark}

Finally, we consider the Schr\"odinger operator $H_{\cG_t}=\D_{\cG_t}+Q$ on the perturbed graph $\cG_t$ obtained from the hexagonal lattice $\cG$ by adding an edge $\be_o=(v_1,v_2)$ with index $\t(\be_o)=t$, $t\in\Z^2$, to its fundamental graph $\cG_*$, see Fig.~\ref{figS2}. Without loss of generality we may assume that the potential $Q$ satisfies
\[\lb{pQhl}
Q_{v_1}=-Q_{v_2}=q>0.
\]

\begin{example}\lb{Exa3} For each $t\in\Z^2$ the perturbed graph $\cG_t$ shown in Fig.~\ref{figS2}\emph{a} and its limit graph $\wt\cG$ (Fig.~\ref{figS3}\emph{b}) are isospectral, and
\[\lb{isE3}
\s(H_{\cG_t})=\s(H_{\wt\cG})=
\big[4-\sqrt{16+q^2}\,,4-q\big]\cup\big[4+q,4+\sqrt{16+q^2}\,\big],
\]
where $H_{\cG_t}=\D_{\cG_t}+Q$ and $H_{\wt\cG}=\D_{\wt\cG}+Q$ are the Schr\"odinger operators on the graphs $\cG_t$ and $\wt\cG$, respectively, with the periodic potential $Q$ defined by \er{pQhl}.
\end{example}

\begin{remark}
These examples also show that perturbations of a periodic graph $\cG$ by adding an edge to its fundamental graph $\cG_*$ essentially change the spectrum of the Schr\"odinger operator $H_{\cG}$ on the unperturbed graph $\cG$ (compare, for example, \er{spHL} and \er{stin}). The spectra of the perturbed problems occur to be closer to the spectrum of the Schr\"odinger operator $H_{\wt\cG}$ on a new periodic graph $\wt\cG$ of a higher dimension (the limit graph) than to the spectrum of the unperturbed problem.
\end{remark}

\subsection{Notions of isospectrality for periodic graphs}\lb{Siso}
There are some notions of isospectrality for periodic graphs. Let $H_1=\D+Q_1$ and $H_2=\D+Q_2$ be two Schr\"odinger operators with periodic potentials $Q_1$ and $Q_2$ on a periodic graph $\cG$.

$\bu$ $H_1$ and $H_2$ are \emph{Floquet isospectral} if
$$
\s(H_1(k))=\s(H_2(k)) \qqq \textrm{for any} \qqq k\in\T^d,
$$
where $H_s(k)$, $s=1,2$, are the corresponding fiber (or Floquet) operators, see \er{Hvt'}, \er{fado}.

$\bu$ $H_1$ and $H_2$ are \emph{Fermi isospectral} if
$$
F_\l(H_1)=F_\l(H_2) \qqq \textrm{for some} \qqq \l\in\R,
$$
where the \emph{Fermi surface} $F_\l(H_s)=\{k\in\T^{d}:\l\in\s(H_s(k))\}$, $s=1,2$.

$\bu$ $H_1$ and $H_2$ are \emph{periodic isospectral} if
$\s(H_1(0))=\s(H_2(0))$, i.e., $H_1$ and $H_2$ have the same \emph{periodic} spectrum.

Floquet isospectrality is the strongest one: if $H_1$ and $H_2$ are Floquet isospectral, then $H_1$ and $H_2$ are Fermi isospectral, periodic isospectral and just isospectral, i.e., $\s(H_1)=\s(H_2)$. The problem of Floquet isospectrality for the discrete Schr\"odinger operators on the lattice $\Z^d$ was studied in \cite{GKT93,Ka89,L23}. Fermi isospectrality on $\Z^d$ was discussed in \cite{L24}.

\begin{figure}[t!]\centering

\includegraphics[width=5cm,height=5.5cm]{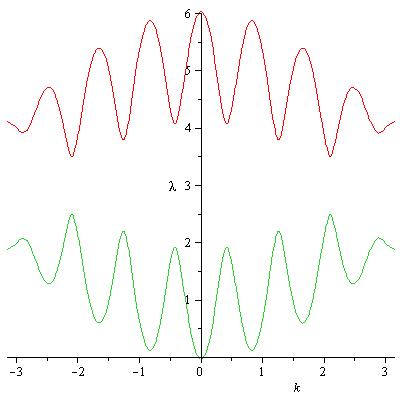}\hspace{20mm}
\includegraphics[width=6cm,height=5.5cm]{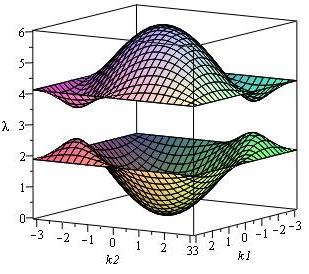}
\caption{\scriptsize The graphs of the band functions (the dispersion relation) of $H_{\dL_t}$ when $t=7$ (left) and the graphs of the band functions of $H_{\wt\cG}$ (right); $q=\frac12$. The dispersion relation for $H_{\dL_t}$ is the restriction of the dispersion relation for $H_{\wt\cG}$ to the line $k_1t-k_2=0$. The ranges of the corresponding band functions coincide.}
\label{fDR}
\end{figure}

Isospectrality of periodic graphs in the sense of coinciding their spectra as sets does not imply Floquet and even Fermi isospectrality, since it only means that the ranges of the corresponding band functions are the same, but the band functions themselves may even depend on the quasimomenta of distinct dimensions. For example, if $t-1\in3\Z$, then the perturbed one-dimensional lattice $\dL_t$ shown in Fig.~\ref{fig5}\emph{b} is isospectral to its limit graph $\wt\cG$, Fig.~\ref{fig5}\emph{c} (see Example~\ref{Ehex}). The graphs of the band functions (the \emph{dispersion relations}) of the Schr\"odinger operators $H_{\dL_t}$ on $\dL_t$ (when $t=7$) and $H_{\wt\cG}$ on $\wt\cG$ are plotted in Fig.~\ref{fDR} (as $q=\frac12$). The band functions of $H_{\dL_t}$ (left) are functions of the one-dimensional quasimomentum $k$. The band functions of $H_{\wt\cG}$ (right) depend on the two-dimensional quasimomentum $k=(k_1,k_2)$. But the ranges of these band functions are the same.

\begin{remark}
\emph{i}) For each $t=(t_1,\ldots,t_d)\in\Z^d$ the perturbed graph $\cG_t$ and its limit graph $\wt\cG$ have the same periodic spectrum (i.e., they are also periodic isospectral). Indeed, the fiber operators $H_{\cG_t}(0)$ and $H_{\wt\cG}(0)$ are the Schr\"odinger operators with the same potential $Q$ on the same fundamental graph $\wt\cG_*$ (see Remark \ref{Rfo0}).

\emph{ii}) The dispersion relation for the Schr\"odinger operator $H_{\cG_t}$ on the perturbed graph $\cG_t$ is just the restriction of the dispersion relation for the Schr\"odinger operator $H_{\wt\cG}$ on the limit graph $\wt\cG$ to the hyperplane $k_1t_1+\ldots+k_dt_d-k_{d+1}=0$, see Proposition \ref{Pnin0}.
\end{remark}

\section{Proof of the main results}\lb{Sec3}
\setcounter{equation}{0}
\subsection{Fiber Schr\"odinger operators for perturbed and limit graphs.} We obtain a relation between the fiber Schr\"odinger operators for the perturbed graph $\cG_t$ ($t\in\Z^d$) and its limit graph $\wt\cG$. Recall that if the perturbed graph $\cG_t$ has dimension $d$, then its limit graph $\wt\cG$ has dimension $d+1$. By construction, the graphs $\cG_t$ and $\wt\cG$ have the same fundamental graph
$$
\wt\cG_*=(\cV_*,\wt\cA_*),\qqq\textrm{where}\qqq
\wt\cA_*=\cA_*\cup\{\be_o,\ol\be_o\},
$$
and $\be_o$ is an edge with index $\t(\be_o)=t$ added to the unperturbed fundamental graph $\cG_*=(\cV_*,\cA_*)$. Here $\cA_*$ and $\wt\cA_*$ are the sets of doubled oriented edges of $\cG_*$ and $\wt\cG_*$, respectively.

Let $\cP:\R^{d+1}\to\R^d$ be the orthogonal projection defined by
\[\lb{cPx}
\cP x=(x_1,\ldots,x_d),\qqq \forall\, x=(x_1,\ldots,x_d,x_{d+1})\in\R^{d+1}.
\]
The adjoint operator $\cP^*:\R^d\to\R^{d+1}$ is the inclusion operator given by
\[\lb{inop}
\cP^*x=(x_1,\ldots,x_d,0)\in\R^{d+1},\qqq \forall\, x=(x_1,\ldots,x_d)\in\R^d.
\]

By construction of the limit graph $\wt\cG$ (see Section \ref{Sec2.1}), the edge indices $\wt\t:\wt\cA_*\to\Z^{d+1}$ for $\wt\cG$ are given by
\[\lb{wtta}
\wt\t(\be)=\left\{
\begin{array}{cl}
\cP^*\t(\be), & \textrm{if $\be\in\cA_*$} \\[2pt]
\mathfrak{e}_{d+1},\; & \textrm{if $\be=\be_o$}\\[2pt]
-\mathfrak{e}_{d+1},\quad & \textrm{if $\be=\ol\be_o$}
\end{array}\right.,\qq\textrm{where}\qq \mathfrak{e}_{d+1}=(0,\ldots,0,1)\in\Z^{d+1},
\]
$\t(\be)\in\Z^d$ is the index of the edge $\be\in\cA_*$ of $\cG_*=(\cV_*,\cA_*)$, and $\cP^*:\R^d\to\R^{d+1}$ is the inclusion operator given by \er{inop}.

\begin{proposition}\lb{Pnin0}
Let $t=(t_1,\ldots,t_d)\in\Z^d$, and let $H_{\cG_t}(k)=\D_{\cG_t}(k)+Q$, $k\in\T^d$, and $H_{\wt\cG}(\wt k)=\D_{\wt\cG}(\wt k)+Q$, $\wt k\in\T^{d+1}$, be the fiber Schr\"odinger operators given by \er{Hvt'}, \er{fado} for the perturbed graph $\cG_t$ and for its limit graph $\wt\cG$, respectively. Then
\[\lb{fLco}
H_{\cG_t}(k)=H_{\wt\cG}(k,t_1k_1+\ldots+t_dk_d),\qqq \forall\,k=(k_1,\ldots,k_d)\in\T^d,
\]
and, consequently, the band functions $\l_{t,j}(\cdot)$ and $\wt\l_j(\cdot)$, $j\in\N_\n$, of the Schr\"odinger operators $H_{\cG_t}$ on $\cG_t$ and $H_{\wt\cG}$ on $\wt\cG$, respectively, satisfy
\[\lb{ljljo}
\l_{t,j}(k)=\wt\l_j(k,t_1k_1+\ldots+t_dk_d),\qqq \forall\,k\in\T^d.
\]
In particular,
\[\lb{incl}
\l_j^-(t)=\min_{\wt k\in\cB_t}\wt\l_j(\wt k)\geq\wt\l_j^-,\qqq
\l_j^+(t)=\max_{\wt k\in\cB_t}\wt\l_j(\wt k)\leq\wt\l_j^+,
\]
where $\l_j^\pm(t)$ and $\wt\l_j^\pm$ are the band edges of $H_{\cG_t}$ and $H_{\wt\cG}$, correspondingly, and
\[\lb{cbt}
\cB_t=\big\{(k_1,\ldots,k_d,k_{d+1})\in\T^{d+1} : t_1k_1+\ldots+t_dk_d- k_{d+1}\in2\pi\Z\big\}.
\]
\end{proposition}

\no \textbf{Proof.}
Due to \er{fado}, the fiber Laplacian $\D_{\wt\cG}(\wt k)$, $\wt k\in\T^{d+1}$, for the limit graph $\wt\cG$ is given by
\[\lb{fiLw}
\big(\D_{\wt\cG}(\wt k)f\big)_v=
\sum_{\be=(v,u)\in\wt\cA_*}\big(f_v-e^{i\lan\wt\t(\be),\wt k\ran}f_u\big),\qqq f\in\ell^2(\cV_*),\qqq v\in \cV_*,
\]
where $\wt\t(\be)\in\Z^{d+1}$ is the index of the edge $\be\in\wt\cA_*$ for the limit graph $\wt\cG$. Using \er{wtta}, we have
\[\lb{tkES}
\lan\wt\t(\be),\wt k\,\ran=\left\{
\begin{array}{cl}
\lan\t(\be),\cP\wt k\ran, & \textrm{if }\be\in\cA_* \\[2pt]
\;k_{d+1}, & \textrm{if } \be=\be_o \\[2pt]
-k_{d+1}, & \textrm{if } \be=\ol\be_o
\end{array}\right.,\qqq \forall\,\wt k=(k_1,\ldots,k_d,k_{d+1})\in\T^{d+1},
\]
where $\cP:\R^{d+1}\to\R^d$ is the orthogonal projection defined by \er{cPx}.

The added edge $\be_o$ of the perturbed graph $\cG_t$ has the index $\t(\be_o)=t$. If $k_{d+1}=t_1k_1+\ldots+t_dk_d$, then the identity \er{tkES} yields
$$
\lan\wt\t(\be),\wt k\ran=\lan\t(\be),\cP\wt k\ran,\qqq\forall\,\be\in\wt\cA_*,
$$
where $\t(\be)\in\Z^{d}$ is the index of the edge $\be\in\wt\cA_*$ for the perturbed graph $\cG_t$. This, \er{fiLw} and \er{fado} give
$$
\big(\D_{\wt\cG}(k_1,\ldots,k_d,t_1k_1+\ldots+t_dk_d)f\big)_v=
\sum_{\be=(v,u)\in\wt\cA_*}\big(f_v-e^{i\lan\t(\be),\cP\wt k\ran}f_u\big)=
\big(\D_{\cG_t}(\cP\wt k)f\big)_v,
$$
which yields \er{fLco}. The identity \er{ljljo} follows from \er{fLco}.

We prove the first formula in \er{incl}. The proof of the second one is similar. Using the definition of the band edges in \er{ban.1H}, the identity \er{ljljo} and $2\pi\Z^{d+1}$-periodicity of the band functions $\wt\l_j(\cdot)$, we obtain
$$
\l_j^-(t)=\min\limits_{k\in\T^d}\l_{t,j}(k)=
\min\limits_{k=(k_1,\ldots,k_d)\in\T^d}\wt\l_j(k,t_1k_1+\ldots+t_dk_d)=
\min\limits_{\wt k\in\cB_t}\wt\l_j(\wt k)\geq\min_{\wt k\in\T^{d+1}}\wt\l_j(\wt k)=\wt\l_j^-,
$$
where the set $\cB_t\ss\T^{d+1}$ is defined by \er{cbt}. \qq $\Box$

\subsection{Asymptotically isospectral graphs.} We prove Theorems \ref{TAig} and \ref{TJig} about asymptotic isospectrality and just isospectrality of the perturbed graph $\cG_t$ and its limit graph~$\wt\cG$.

\medskip

\no \textbf{Proof of Theorem \ref{TAig}.} Let $j\in\N_\n$. We prove \er{as00} for $\wt\l_j^-$. The proof for $\wt\l_j^+$ is similar. Due to \er{incl}, \er{cbt}, in order to find the band edge $\l_j^-(t)$ we need to minimize the band function $\wt\l_j(\cdot)$ of $H_{\wt\cG}$ over the set of the hyperplanes
\[\lb{hpLn}
L_n: t_1k_1+\ldots+t_dk_d-k_{d+1}=2\pi n, \qqq n\in\Z,
\]
in $\T^{d+1}=(-\pi,\pi]^{d+1}$.

Let $k_o=(k_1^o,\ldots,k_d^o,k_{d+1}^o)\in K_j^-$, where $K_j^-$ is defined in \er{Kjpm}. The distance $\r(n)$ from the point $k_o$ to the hyperplane $L_n$ is equal to
\[\lb{kokh0}
\r(n)=\frac{|t_1k_1^o+\ldots+t_dk_d^o-k_{d+1}^o-2\pi n|}{\sqrt{|t|^2+1}}=\frac{2\pi|\vr_t(k_o)-n|}{\sqrt{|t|^2+1}}\,,
\]
where $\vr_t(\cdot)$ is given by \er{fual}. Let $\wh n:=\mn(\vr_t(k_o))$ be the nearest integer to $\vr_t(k_o)$. Then, among the hyperplanes \er{hpLn},  $L_{\wh n}$ is closest to the extreme point $k_o$. Let $k^\perp$ be the orthogonal projection of $k_o$ onto the hyperplane $L_{\wh n}$ (see Fig. \ref{fig0}). Then, using \er{kokh0}, we obtain
\[\lb{kokh}
|k_o-k^\perp|=\r(\wh n)=\frac{2\pi|\vr(k_o)-\mn(\vr_t(k_o))|}{\sqrt{|t|^2+1}}=
\frac{2\pi|\a(\vr_t(k_o))|}{\sqrt{|t|^2+1}}=\textstyle O\big(\frac1{|t|}\big), \qq \textrm{as}\qq |t|\to\iy,
\]
where $\a(\vr_t(k_o))\in\big(-\frac12,\frac12\big]$ is defined by \er{alox}.
Since $k^\perp\in L_{\wh n}$ and due to the first formula in \er{incl}, we have
$$
\wt\l_j^-\leq\l_j^-(t)=\min_{\wt k\in\cB_t}\wt\l_j(\wt k)\leq\wt\l_j(k^\perp).
$$
This, \er{kokh} and continuity of the band function $\wt\l_j(\cdot)$ yield
$$
\wt\l_j^-\leq\lim_{|t|\to\iy}\l_j^-(t)\leq\lim_{|t|\to\iy}\wt\l_j(k^\perp)
=\wt\l_j(k_o)=\wt\l_j^-.
$$
Thus, $\lim\limits_{|t|\to\iy}\l_j^-(t)=\wt\l_j^-$. \qq $\Box$

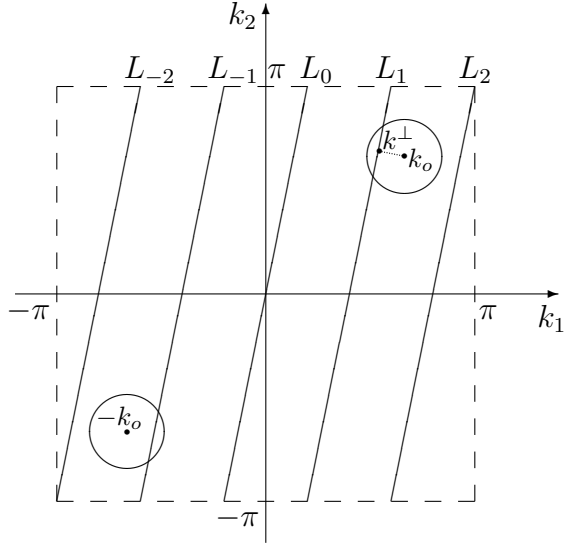
\begin{figure}[t!]\centering
\unitlength 1.1mm 
\linethickness{0.4pt}

\begin{picture}(65,65)
\put(13.4,13.4){\circle*{0.7}}
\put(13.4,13.4){\circle{9}}

\put(46.6,46.6){\circle*{0.7}}
\put(46.6,46.6){\circle{9}}
\put(43.6,47.2){\circle*{0.7}}

\put(46.6,46.6){\line(-5,1){1.0}}
\bezier{10}(43.6,47.2)(45.1,46.9)(46.6,46.6)

\put(9.7,14.2){\footnotesize$-k_o$}
\put(47.0,45.3){\footnotesize$k_o$}
\put(43.9,47.5){\footnotesize$k^{\perp}$}

\put(0,30){\vector(1,0){65.00}}
\put(30,0){\vector(0,1){65.00}}

\put(5,5){\line(1,5){10.0}}
\put(15,5){\line(1,5){10.0}}
\put(25,5){\line(1,5){10.0}}
\put(35,5){\line(1,5){10.0}}
\put(45,5){\line(1,5){10.0}}

\multiput(5,5)(4,0){13}{\line(1,0){2}}
\multiput(5,55)(4,0){13}{\line(1,0){2}}

\multiput(5,5)(0,4){13}{\line(0,1){2}}
\multiput(55,5)(0,4){13}{\line(0,1){2}}

\put(55.5,27){$\pi$}
\put(-1,27){$-\pi$}
\put(30.2,56.0){$\pi$}
\put(24,2){$-\pi$}

\put(13,56){$L_{-2}$}
\put(23,56){$L_{-1}$}
\put(34,56){$L_0$}
\put(43,56){$L_1$}
\put(53,56){$L_2$}
\put(62.5,26){$k_1$}
\put(25.5,62.5){$k_2$}
\end{picture}
\caption{\scriptsize Global extreme points $\pm k_o$ of a band function $\wt\l_j(k)$, $k=(k_1,k_2)\in\T^2=(-\pi,\pi]^2$, and their small neighborhoods; the lines $L_n$ are given by the equations $tk_1-k_2=2\pi n$, $n\in\Z$ (when $t=5$). The line $L_1$ is the line in $\cB_t$ closest to the extreme point $k_o$, and $k^\perp$ is the orthogonal projection of $k_o$ onto $L_1$; $\big|k_o-k^\perp\big|=O(\frac1{|t|})$ as $|t|\to\iy$.}
\label{fig0}
\end{figure}

\medskip

\no \textbf{Proof of Theorem \ref{TJig}.} We prove this theorem for the band edge $\l_j^-(t)$. The proof for $\l_j^+(t)$ is similar.  Let the condition \er{soeq} be fulfilled for some $k_o\in K_j^-$. Then $k_o\in\cB_t$, where $\cB_t$ is defined by \er{cbt}, and, due to \er{incl},
$$
\wt\l_j^-\leq\l_j^-(t)=\min_{\wt k\in\cB_t}\wt\l_j(\wt k)\leq\wt\l_j(k_o)=\wt\l_j^-,
$$
where $\wt\l_j(\cdot)$ is the band function of the Schr\"odinger operators $H_{\wt\cG}$. Thus, $\l_j^-(t)=\wt\l_j^-$.

Conversely, let $\l_j^-(t)=\wt\l_j^-$. Then, due to \er{incl}, for some $k_o\in\cB_t$ we have
$$
\l_j^-(t)=\min_{\wt k\in\cB_t}\wt\l_j(\wt k)=\wt\l_j(k_o)=\wt\l_j^-.
$$
Thus, $k_o\in K_j^-$ and $k_o$ satisfies the condition \er{soeq} (since $k_o\in\cB_t$).

Finally, let $0\in K_j^-$. For $k_o=0$ the condition \er{soeq} holds. Then $\l_j^-(t)=\wt\l_j^-$.\qq $\Box$

\subsection{Asymptotics of the band edges.} We prove Theorem \ref{TAsy} about asymptotics of the band edges $\l_j^\pm(t)$ of the Schr\"odinger operator $H_{\cG_t}$ on the perturbed graph $\cG_t$ as $|t|\to\iy$.

\medskip

\no \textbf{Proof of Theorem \ref{TAsy}.} We prove \er{aslk} for $\l_j^-(t)$. The proof for $\l_j^+(t)$ is similar. Due to \er{incl},
\[\lb{Tla-}
\l_j^-(t)=\min_{k\in\cB_t}\wt\l_j(k),
\]
where $\cB_t$ is defined by \er{cbt}. Thus, in order to find the band edge $\l_j^-(t)$ we need to minimize $\wt\l_j(\cdot)$ over the set of the hyperplanes
\[\lb{cong0}
L_n: \lan \gt,k\ran=2\pi n, \qq n\in\Z,\qqq \textrm{where}\qqq \gt=(t_1,\ldots,t_d,-1),
\]
in $\T^{d+1}=(-\pi,\pi]^{d+1}$. For the single (up to evenness) global minimum point $k_o\in\T^{d+1}$ of the function $\wt\l_j(k)$ we find the global minimum of $\wt\l_j(k)$ with the constraint \er{cong0} in a small neighborhood of $k_o$, see Fig. \ref{fig0}.

Assumption A1 means that the Hessian of the band function $\wt\l_j$ at $k_o$ is positive definite. The Taylor expansion of the function $\wt\l_j(k)$ about the point $k_o$ is given by
\[\lb{wtlk0}
\begin{array}{l}
\textstyle\wt\l_j(k)=\wt\l_j^-+\frac12\lan k-k_o,\bH(k-k_o)\ran+g_1(k-k_o),\qqq\textrm{where} \\[6pt]
\wt\l_j^-=\wt\l_j(k_o),\qqq \bH=\mathrm{Hess}\,\wt\l_j(k_o),\qqq g_1(k-k_o)=O\big(|k-k_o|^3\big).
\end{array}
\]

We make the non-degenerate change of variables $k\mapsto y$, where
$$
y=\bH^{1/2}\,(k-k_o),
$$
and $\bH^{1/2}$ is the positive definite square root of $\bH$. In these new variables $y$ the expansion \er{wtlk0} and the constraint \er{cong0} have the form
\[\lb{wtlk}
\textstyle\wt\l_j(k(y))=\wt\l_j^-+\frac12|y|^2
+g_2(y),\qqq g_2(y)=g_1(\bH^{-1/2}y)=O\big(|y|^3\big),
\]
\[\lb{hupy}
\lan \bH^{-1/2}\gt,y\ran=2\pi n-\lan \gt,k_o\ran, \qqq n\in\Z.
\]

Using the Lagrange multiplier method one can show that for a linear equation $\lan a,y\ran=b$, where $a,y\in\R^n$, $b\in\R$, the solution $y_o$ of minimum Euclidean norm $|y_o|$ is given by $y_o=\frac{b}{|a|^2}\,a$. Then for each $n\in\Z$ the solution $y_o$ of the linear equation \er{hupy} with minimum $|y_o|$ has the form
$$
y_o=y_o(n)=\frac{2\pi n-\lan \gt,k_o\ran}{|\bH^{-1/2}\gt|^2}\,\bH^{-1/2}\gt,
$$
or, using \er{cong0},
\[\lb{noys}
y_o(n)=\frac{2\pi(n-\vr_t(k_o))}{|\bH^{-1/2}\gt|^2}\,\bH^{-1/2}\gt,
\]
where $\vr_t(\cdot)$ is given by \er{fual}. We find $\min\limits_{n\in\Z}|y_o(n)|$. Using \er{noys}, we have
\[\lb{nyo}
|y_o(n)|=\dfrac{2\pi|\vr_t(k_o)-n|}{|\bH^{-1/2}\gt|}\,.
\]
Since the minimum of $|\vr_t(k_o)-n|$ is achieved for $n=\wh n:=\mn(\vr_t(k_o))$, where $\mn(x)$ is the nearest integer to $x$, then \er{nyo} yields
\[\lb{omal}
\min_{n\in\Z}|y_o(n)|=|y_o(\wh n)|
=\frac{2\pi|\vr_t(k_o)-\wh n|}{|\bH^{-1/2}\gt|}=
\frac{2\pi|\a(\vr_t(k_o))|}{|\bH^{-1/2}\gt|}=\textstyle O(\frac1{|t|}),
\]
where $\a(x)\in\big(-\frac12,\frac12\big]$ is defined by \er{alox}. Let $\wh y_o=y_o(\wh n)$. Then, using \er{Tla-}, \er{wtlk} and \er{omal}, we obtain
$$
\textstyle\l_j^-(t)=\wt\l_j(k(\wh y_o))=\wt\l_j^-+\frac12\,|\wh y_o|^2
+g_2(\wh y_o)=\wt\l_j^-+\dfrac{2\pi^2\a^2(\vr_t(k_o))}{|\bH^{-1/2}\gt|^2}+g_2(\wh y_o),
$$
where $g_2(\wh y_o)=O(|\wh y_o|^3)$, or, using \er{omal}, $g_2(\wh y_o)=O(\frac1{|t|^3})$. Thus, the asymptotics \er{aslk} is proved.

If all components of $k_o$ are equal to either 0 or $\pi$, then, due to $2\pi\Z^{d+1}$-periodicity and evenness of $\wt\l_j$, we have
$$
\wt\l_j(k+k_o)=\wt\l_j(-k-k_o)=\wt\l_j(-k+k_o),
$$
since $2k_o\in2\pi\Z^{d+1}$. This implies that the Taylor series \er{wtlk0} of $\wt\l_j$ at $k_o$ has only even degree terms. Then, $g_1(k-k_o)=O\big(|k-k_o|^4\big)$ and $g_2(\wh y_o)=O(\frac1{|t|^4})$. \qq $\Box$

\section{Examples}
\setcounter{equation}{0}
\lb{Sec4}
\subsection{Perturbations of the one-dimensional lattice}\lb{Ss4.2}
We consider the Schr\"odinger operator $H_\dL=\D_\dL+Q$ on the one-dimensional lattice $\dL$ (Fig.~\ref{fig5}\emph{a}) with a 2-periodic potential $Q$ defined by $Q_0=-Q_1=q>0$. The fundamental graph $\dL_*=\dL/2\Z$ of the lattice $\dL$ consists of two vertices $v_0$ and $v_1$ with degree $\vk_{v_0}=\vk_{v_1}=2$ and two edges $\be_1,\be_2$ connecting these vertices, see Fig.~\ref{fig5}\emph{a} (right). The indices of the edges are given by
$$
\t(\be_1)=0, \qqq \t(\be_2)=1.
$$
The fiber Schr\"odinger operator $H_\dL(k)$ defined by \er{Hvt'}, \er{fado} for the graph $\dL$ has the form
$$
H_\dL(k)=\left(
\begin{array}{cc}
  2+q & -1-e^{ik}\\
  -1-e^{-ik} & 2-q
\end{array}\right), \qqq k\in\T=(-\pi,\pi].
$$
The eigenvalues of $H_\dL(k)$ are given by
$$
\textstyle\l_j(k)=2+(-1)^j\sqrt{4\cos^2\frac k2+q^2}\,,\qqq j=1,2.
$$
Thus,
$$
\begin{array}{ll}
\textstyle \l_1^-=\min\limits_{k\in\T}\l_1(k)=\l_1(0)=2-\sqrt{4+q^2}, \qq & \l_1^+=\max\limits_{k\in\T}\l_1(k)=\l_1(\pi)=2-q,\\[12pt]
\textstyle \l_2^-=\min\limits_{k\in\T}\l_2(k)=\l_2(\pi)=2+q, \qq & \l_2^+=\max\limits_{k\in\T}\l_2(k)=\l_2(0)=2+\sqrt{4+q^2},
\end{array}
$$
and the spectrum of the Schr\"odinger operator $H_\dL$ has the form \er{spHL}.

\medskip

Let $t\in\Z$, and let $\dL_t$ be the perturbed graph obtained from the lattice $\dL$ by adding edges $(2n,2n+1+2t)$, $n\in\Z$, see Fig.~\ref{fig5}\emph{b}. The perturbed graph $\dL_t$ has the period 2. The fundamental graph $\wt\dL_*$ of $\dL_t$ consists of the edges $\be_1,\be_2$ and the added edge $\be_o$ connecting the vertices $v_0$ and $v_1$, see Fig.~\ref{fig5}\emph{b} (right). The index of the added edge $\be_o$ is $\t(\be_o)=t$. The limit graph of $\dL_t$ has the following construction (see Section \ref{Sec2.1}). Let $\ga_2=(0,1)$. We take an infinite number of copies
$$
\dL+n\ga_2\ss\mathbb{R}^{2},\qqq \forall\,n\in\Z,
$$
of the one-dimensional lattice $\dL$ (shown in black color in Fig. \ref{fig5}\emph{c}) and connect them by edges
$$
\be_o+\ga,\qq \forall\,\ga\in\wt\G,\qqq \textrm{where}\qqq
\be_o=(v_0,v_1+\ga_2),
$$
(red color in Fig. \ref{fig5}\emph{c}). Here $\wt\G$ is the lattice in $\R^{2}$ with the basis $\ga_1=(2,0)$, $\ga_2=(0,1)$. The obtained graph $\wt\cG$ is the limit graph for the perturbed graph $\dL_t$.

\medskip

First we describe the spectrum of the Schr\"odinger operator $H_{\wt\cG}=\D_{\wt\cG}+Q$ with a $\wt\G$-periodic potential $Q$ on the limit graph~$\wt\cG$.

\begin{lemma}\lb{LLG2}
Let $\wt\cG$ be the graph shown in Fig.~\ref{fig5}c. Then the spectrum of the Schr\"odinger operator $H_{\wt\cG}=\D_{\wt\cG}+Q$ with the periodic potential $Q$ defined by $Q_{v_0}=-Q_{v_1}=q>0$ on $\wt\cG$ is given by
\[\lb{sDbG}
\begin{aligned}
&\s(H_{\wt\cG})=[\,\wt\l_1^-,\wt\l_1^+]\cup[\,\wt\l_2^-,\wt\l_2^+],\\
& \begin{array}{ll}
\wt\l_1^-=\wt\l_1(K_1^-)=3-\sqrt{9+q^2},\qq &  \wt\l_1^+=\wt\l_1(K_1^+)=3-q,\\[6pt]
\wt\l_2^-=\wt\l_2(K_2^-)=3+q, \qq & \wt\l_2^+=\wt\l_2(K_2^+)=3+\sqrt{9+q^2},
\end{array}
\end{aligned}
\]
where the level sets
\[\lb{spDoK}
\textstyle K_1^-=K_2^+=\big\{(0,0)\big\}, \qqq
K_1^+=K_2^-=\big\{\big(\pm\frac{2\pi}3\,,\pm\frac{2\pi}3\big)\big\}.
\]
\end{lemma}
\no \textbf{Proof.} The fundamental graph $\wt\cG_*$ of $\wt\cG$ consists of two vertices $v_0$ and $v_1$ with degree $\vk_{v_0}=\vk_{v_1}=3$ and three edges $\be_1,\be_2,\be_o$ connecting these vertices, see Fig.~\ref{fig5}\emph{c} (right). The indices of the edges are given by
$$
\wt\t(\be_1)=(0,0), \qqq \wt\t(\be_2)=(1,0),\qqq \wt\t(\be_o)=(0,1).
$$
The fiber Schr\"odinger operator $H_{\wt\cG}(k)$ defined by \er{Hvt'}, \er{fado} for the graph $\wt\cG$ has the form
$$
H_{\wt\cG}(k)=\left(
\begin{array}{cc}
  3+q & -1-e^{ik_1}-e^{-ik_2} \\
  -1-e^{-ik_1}-e^{ik_2} & 3-q
\end{array}\right), \qq k=(k_1,k_2)\in\T^2=(-\pi,\pi]^2.
$$
The eigenvalues of $H_{\wt\cG}(k)$ are given by
\[\lb{evEG}
\begin{aligned}
&\wt\l_j(k)=3+(-1)^j\sqrt{f(k)+q^2},\qq j=1,2,\qq \textrm{where}\\
&f(k)=|1+e^{ik_1}+e^{-ik_2}|^2=3+2\cos k_1+2\cos k_2+2\cos(k_1+k_2).
\end{aligned}
\]
We have
$$
\textstyle\min\limits_{k\in\T^2}f(k)=f\big(\pm\frac{2\pi}3\,,\pm\frac{2\pi}3\big)=0,
\qqq
\max\limits_{k\in\T^2}f(k)=f(0,0)=9.
$$
Then, using the definition of the band edges $\wt\l_j^\pm$ in \er{debe},
we obtain that the spectrum of $H_{\wt\cG}$ has the form \er{sDbG}, \er{spDoK}. \qq $\Box$

\medskip

Now we prove Example \ref{Ehex} about the spectrum of the Schr\"odinger operator on the perturbed graph $\dL_t$ shown in Fig.~\ref{fig5}\emph{b}.

\medskip

\no \textbf{Proof of Example \ref{Ehex}.} The limit graph for the perturbed graph $\dL_t$ is the graph $\wt\cG$ shown in Fig.~\ref{fig5}\emph{c}. Then,
by Theorem \ref{TAig}, $\dL_t$ is asymptotically isospectral to $\wt\cG$, i.e.,
$$
\s(H_{\dL_t})=[\l_1^-(t),\l_1^+(t)]\cup[\l_2^-(t),\l_2^+(t)],\qqq
\lim\limits_{|t|\to\iy}\l_j^\pm(t)=\wt\l_j^\pm,\qq j=1,2,
$$
where $\wt\l_j^\pm=\wt\l_j(K_j^\pm)$ are the band edges for the Schr\"odinger operator $H_{\wt\cG}=\D_{\wt\cG}+Q$ on $\wt\cG$, see \er{sDbG}, \er{spDoK}. Since $(0,0)\in K_1^-=K_2^+$, Theorem \ref{TJig} yields
\[\lb{pgep}
\l_1^-(t)=\wt\l_1^-=3-\sqrt{9+q^2},\qqq \l_2^+(t)=\wt\l_2^+=3+\sqrt{9+q^2}.
\]

\emph{i}) Let $t-1\in3\mathbb{Z}$. Then for the points $\big(\pm\frac{2\pi}3\,,\pm\frac{2\pi}3\big)\in K_1^+=K_2^-$ the condition \er{soeq} is fulfilled. Thus, due to Theorem \ref{TJig} and the identities \er{sDbG},
$$
\textstyle\l_1^+(t)=\wt\l_1^+=3-q,\qqq
\l_2^-(t)=\wt\l_2^-=3+q.
$$
This and \er{pgep} yield $\s(H_{\dL_t})=\s(H_{\wt\cG})$.

\emph{ii}) Let $t-1\not\in3\mathbb{Z}$, i.e., $t=3m$ or $t=3m+2$ for some $m\in\Z$. We obtain the asymptotics of the band edge $\l_1^+(t)$ as $|t|\to\iy$. The proof for $\l_2^-(t)$ is similar. The point $k_o=\big(\frac{2\pi}3\,,\frac{2\pi}3\big)$ is the only (up to evenness) maximum point of the band function $\wt\l_1(k)$ in $\T^2$, see \er{sDbG}, \er{spDoK}. Using the identity \er{evEG} for $\wt\l_1(k)$, we obtain
$$
\frac{\pa^2\wt\l_1}{\pa k_j^2}\,(k_o)=-\frac1q\,,\qq j=1,2;
\qqq\frac{\pa^2\wt\l_1}{\pa k_1\pa k_2}\,(k_o)=-\frac1{2q}\,.
$$
Then the matrix $\bH=-\mathrm{Hess}\,\wt\l_1(k_o)$ and its inverse are given by
$$
\bH=\frac1{2q}\left(
\begin{array}{cc}
  2 & 1 \\
  1 & 2
\end{array}\right), \qqq \bH^{-1}=\frac{2q}3\left(
\begin{array}{rr}
  2 & -1 \\
  -1 & 2
\end{array}\right).
$$
Thus, the band function $\wt\l_1(k)$ of $H_{\wt\cG}$ satisfies Assumption A and we can use the asymptotics \er{aslk}. For $\gt=(t,-1)$ we obtain
\[\lb{coef1}
\textstyle|\bH^{-1/2}\gt|^2=\lan\gt,\bH^{-1}\gt\ran=\frac{4q}3\,(t^2+t+1).
\]
Using \er{fual}, we get
$$
\textstyle\vr_t(k_o)=\frac1{2\pi}(t\,\frac{2\pi}3-\frac{2\pi}3)=\frac{t-1}3\,.
$$
Since $t=3m$ or $t=3m+2$ for some $m\in\Z$, the nearest integer to $\vr_t(k_o)$ is $\mn\big(\vr_t(k_o)\big)=m$, and
\[\lb{slal2}
\textstyle \a\big(\vr_t(k_o)\big)=\vr_t(k_o)-\mn\big(\vr_t(k_o)\big)=\left\{
\begin{array}{rl}
-\frac13, & \textrm{if $t=3m$} \\[2pt]
\frac13, & \textrm{if $t=3m+2$}
\end{array}\right.\,,
\]
see \er{alox}. Substituting \er{slal2}, \er{coef1} and $\wt\l_1^+=3-q$ into \er{aslk}, we obtain
$$
\textstyle\l_1^+(t)=3-q-\frac{\pi^2}{6q(t^2+t+1)}+O(\frac1{|t|^3})=
3-q-\frac{\pi^2}{6q}\cdot\frac1{t^2}+O(\frac1{|t|^3}),\qq\textrm{as}\qq |t|\to\iy.\qqq \Box
$$

\subsection{Perturbations of the $d$-dimensional lattice} The $d$-dimensional lattice $\dL^d$ is a $\Z^d$-periodic graph, for $d=2,3$ see  Fig.~\ref{fig2}\emph{a,c}. The fundamental cell $\Omega=[0,1)^d$ is shaded in the figure. The fundamental graph $\dL^d_*$ of the lattice $\dL^d$ consists of one vertex $v$ with degree $\vk_v=2d$ and $d$ loop edges $\be_1,\ldots,\be_d$ at this vertex $v$ with indices
$$
\t(\be_1)=\ga_1,\qq \ldots,\qq \t(\be_d)=\ga_d,
$$
where $\ga_1,\ldots,\ga_d$ is the standard basis of $\Z^d$.

The fiber Laplacian $\D_{\dL^d}(k)$ defined by \er{fado} for the lattice $\dL^d$ has the form
\[\lb{bfDd}
\wt\l_1(k)\equiv\D_{\dL^d}(k)=2d-2\cos k_1-\ldots-2\cos k_d, \qq k=(k_1,\ldots,k_d)\in\T^d=(-\pi,\pi]^d.
\]
Thus,
$$
\textstyle \wt\l_1^-=\min\limits_{k\in\T^d}\wt\l_1(k)=\wt\l_1(K_1^-)=0,  \qqq \wt\l_1^+=\max\limits_{k\in\T^d}\wt\l_1(k)=\wt\l_1(K_1^+)=4d,
$$
where the level sets $K_1^\pm\ss\T^d$ consist of the single point:
$$
K_1^-=\big\{(0,\ldots,0)\big\}, \qqq
K_1^+=\big\{\pi\1_d\big\},\qqq \1_d=(1,\ldots,1)\in\Z^d.
$$
The spectrum of the Laplacian $\D_{\dL^d}$ on the $d$-dimensional lattice $\dL^d$ is given by
\[\lb{sDLd}
\s(\D_{\dL^d})=[\wt\l_1^-,\wt\l_1^+]=[0,4d].
\]

Now we prove Example \ref{ExSl} about the Laplacian spectrum under the perturbations of the lattice $\dL^d$ by adding an edge to its fundamental graph $\dL_*^d$.

\medskip

\no \textbf{Proof of Example \ref{ExSl}.} The limit graph for the perturbed $d$-dimensional lattice $\dL^d_t$ is the $(d+1)$-dimensional lattice $\dL^{d+1}$. Indeed, stacking together infinitely many copies of the $d$-dimensional lattice $\dL^d\ss\R^d$ along the vector $\ga_{d+1}=(0,\ldots,0,1)\in\Z^{d+1}$ and connecting them by edges $(v,v+\ga_{d+1})$, $\forall v\in\Z^{d+1}$, we obtain the lattice $\dL^{d+1}$ (for $d=2$ see Fig.~\ref{fig2}\emph{c}).

The spectrum of the Laplacian $\D_{\dL^{d+1}}$ on $\dL^{d+1}$ has the form
$$
\s(\D_{\dL^{d+1}})=[\wt\l_1^-,\wt\l_1^+]=[0,4(d+1)],
$$
see \er{sDLd}. By Theorem \ref{TAig}, $\dL^d_t$ is asymptotically isospectral to $\dL^{d+1}$, i.e.,
$$
\s(\D_{\dL^d_t})=[\l_1^-(t),\l_1^+(t)],\qqq
\lim\limits_{|t|\to\iy}\l_1^-(t)=0,\qqq
\lim\limits_{|t|\to\iy}\l_1^+(t)=4(d+1).
$$
Since $0\in K_1^-$, then Theorem \ref{TJig} yields $\l_1^-(t)=\wt\l_1^-=0$. Thus, \er{sDLt} is proved.

\emph{i}) Let $t_1+\ldots+t_d$ be odd. Then for $k_o=\pi\1_{d+1}\in K_1^+$ the condition \er{soeq} is fulfilled. Thus, due to Theorem \ref{TJig}, $\l_1^+(t)=\wt\l_1^+=4(d+1)$ and, consequently, $\s(\D_{\dL^d_t})=\s(\D_{\dL^{d+1}})$.

\emph{ii}) Let $t_1+\ldots+t_d$ be even. The single band function of $\D_{\dL^{d+1}}$
$$
\wt\l_1(k)=2(d+1)-2\cos k_1-\ldots-2\cos k_{d+1}, \qqq k=(k_1,\ldots,k_{d+1})\in\T^{d+1},
$$
see \er{bfDd}. Then for the singe maximum point $k_o=\pi\1_{d+1}\in K_1^+$ of $\wt\l_1(\cdot)$ we obtain
$$
\frac{\pa^2 \wt\l_1}{\pa k_i^2}\,(k_o)=-2,\qqq \frac{\pa^2\wt\l_1}{\pa k_i\pa k_j}\equiv0,\qqq
i,j\in\N_{d+1}, \qqq i\neq j.
$$
The matrix $\bH=-\mathrm{Hess}\,\wt\l_1(k_o)$ and its inverse are given by
$$
\textstyle\bH=2I_{d+1},\qqq \bH^{-1}=\frac12I_{d+1},
$$
where $I_{d+1}$ is the identity matrix of size $d+1$. Thus, the band function $\wt\l_1(k)$ of $\D_{\dL^{d+1}}$ satisfies Assumption A and we can use the asymptotics \er{aslk}. For $\gt=(t_1,\ldots,t_d,-1)$ we obtain
\[\lb{coef}
\textstyle|\bH^{-1/2}\gt|^2=\lan\gt,\bH^{-1}\gt\ran=\frac12\,|\gt|^2=\frac12(|t|^2+1).
\]
Using \er{fual}, we get
$$
\textstyle \vr_t(k_o)=\vr_t(\pi\1_{d+1})=\frac{t_1+\ldots+t_d}2-\frac12\,.
$$
Since $t_1+\ldots+t_d$ is even, the nearest integer to $\vr_t(k_o)$ is $\mn\big(\vr_t(k_o)\big)=\frac{t_1+\ldots+t_d}2-1$, and
\[\lb{slal}
\textstyle \a\big(\vr_t(k_o)\big)=\vr_t(k_o)-\mn\big(\vr_t(k_o)\big)=\frac12\,,
\]
see \er{alox}. Substituting \er{slal}, \er{coef} and $\wt\l_1^+=4(d+1)$ into \er{aslk} and using that each component of $k_o$ is $\pi$, we get
$$
\textstyle\l_1^+(t)=4(d+1)-\frac{\pi^2}{|t|^2+1}+O\big(\frac1{|t|^4}\big)
=4(d+1)-\frac{\pi^2}{|t|^2}+O\big(\frac1{|t|^4}\big),
\qqq \textrm{as}\qqq |t|\to\iy.
$$
Thus, the asymptotics \er{ase1} is proved. \qq $\Box$

\subsection{Perturbations of the hexagonal lattice}
Let $H_\cG=\D_\cG+Q$ be the Schr\"odinger operator on the hexagonal lattice $\cG$ (see Fig. \ref{figS1}). The fundamental graph $\cG_*$ of $\cG$ consists of two vertices $v_1$ and $v_2$ and three multiple edges connecting these vertices. Without loss of generality we may assume that
\[\lb{poQs}
Q_{v_1}=-Q_{v_2}=q>0.
\]
Since the hexagonal lattice $\cG$ is isomorphic to the graph shown in Fig.~\ref{fig5}\emph{c}, due to Lemma \ref{LLG2}, the spectrum of the Schr\"odinger operator $H_\cG$ has the form
$$
\s(H_\cG)=\big[3-\sqrt{9+q^2}\,,3-q\big]\cup\big[3+q,3+\sqrt{9+q^2}\,\big].
$$

Let $\cG_t$ be the perturbed graph obtained from the hexagonal lattice $\cG$ by adding an edge $\be_o=(v_1,v_2)$ with index $\t(\be_o)=t$ ($t\in\Z^2$) to its fundamental graph $\cG_*$, see Fig.~\ref{figS2}. By construction (see Section \ref{Sec2.1}), the graph $\wt\cG$ shown in Fig.~\ref{figS3}\emph{b} is the limit graph for the perturbed hexagonal lattice $\cG_t$. First we describe the spectrum of the Schr\"odinger operator $H_{\wt\cG}=\D_{\wt\cG}+Q$ on the limit graph~$\wt\cG$.

\begin{lemma}\lb{Ldeb}
Let $\wt\cG$ be the graph shown in Fig.~\ref{figS3}b. Then the spectrum of the Schr\"odinger operator $H_{\wt\cG}=\D_{\wt\cG}+Q$ with the potential $Q$ defined by \er{poQs} on $\wt\cG$ is given by
\[\lb{swbG}
\begin{aligned}
&\s(H_{\wt\cG})=[\,\wt\l_1^-,\wt\l_1^+]\cup[\,\wt\l_2^-,\wt\l_2^+],\\
& \begin{array}{ll}
\wt\l_1^-=\wt\l_1(K_1^-)=4-\sqrt{16+q^2},\qq &  \wt\l_1^+=\wt\l_1(K_1^+)=4-q,\\[6pt]
\wt\l_2^-=\wt\l_2(K_2^-)=4+q, \qq & \wt\l_2^+=\wt\l_2(K_2^+)=4+\sqrt{16+q^2},
\end{array}
\end{aligned}
\]
where the level sets $K_1^-=K_2^+=\big\{(0,0,0)\big\}$,
\[\lb{sbG3}
\textstyle  K_1^+=K_2^-=\big\{(k_1,k_2,k_3)\in\T^3 : 1+e^{ik_1}+e^{ik_2}+e^{ik_3}=0\big\}.
\]
\end{lemma}

\begin{remark}
The set $K_1^+=K_2^-$ defined by \er{sbG3} is an algebraic curve of co-dimension 2, and the band edges $\wt\l_1^+$ and $\wt\l_2^-$ are degenerate. This example of periodic graphs where the non-degeneracy assumption fails was constructed in \cite{BK20}.
\end{remark}

\no \textbf{Proof.} The fundamental graph $\wt\cG_*$ of $\wt\cG$ consists of two vertices $v_1$ and $v_2$ with degree $\vk_{v_1}=\vk_{v_2}=4$ and four edges $\be_1,\be_2,\be_3,\be_o$ connecting these vertices, see Fig.~\ref{figS2}\emph{b}. The indices of the edges are given by
$$
\wt\t(\be_1)=(0,0,0), \qqq \wt\t(\be_2)=(1,0,0), \qqq \wt\t(\be_3)=(0,1,0),\qqq \wt\t(\be_o)=(0,0,1).
$$
The fiber Schr\"odinger operator $H_{\wt\cG}(k)$, $k=(k_1,k_2,k_3)\in\T^3$, defined by \er{Hvt'}, \er{fado} for the graph $\wt\cG$ has the form
$$
H_{\wt\cG}(k)=\left(
\begin{array}{cc}
  4+q & -1-e^{-ik_1}-e^{-ik_2}-e^{-ik_3} \\
  -1-e^{ik_1}-e^{ik_2}-e^{ik_3} & 4-q
\end{array}\right).
$$
The eigenvalues of $H_{\wt\cG}(k)$ are given by
$$
\begin{aligned}
&\wt\l_j(k)=4+(-1)^j\sqrt{f(k)+q^2},\qqq j=1,2,\qq  \\
&\textrm{where} \qqq f(k)=|1+e^{ik_1}+e^{ik_2}+e^{ik_3}|^2.
\end{aligned}
$$
We have
$$
\textstyle\min\limits_{k\in\T^3}f(k)=f\big(K_1^+\big)=0,\qqq
\max\limits_{k\in\T^3}f(k)=f(0,0,0)=16,
$$
where $K_1^+$ is defined by \er{sbG3}. Then, using the definition of the band edges $\wt\l_j^\pm$ in \er{debe}, we obtain that the spectrum of $H_{\wt\cG}$ has the form \er{swbG}, \er{sbG3}. \qq $\Box$

\medskip

Finally, we prove Example \ref{Exa3} about isospectrality of the graphs $\cG_t$ and $\wt\cG$ shown in Fig.~\ref{figS2}\emph{a} and Fig.~\ref{figS3}\emph{b}, respectively.

\medskip

\no \textbf{Proof of Example \ref{Exa3}.} Let $t=(t_1,t_2)\in\Z^2$, and let $k_o\in\T^3$ be defined by
$$
k_o=\left\{
\begin{array}{rl}
(\pi,0,\pi), & \textrm{if $t_1$ is odd,} \\[2pt]
(0,\pi,\pi), & \textrm{if $t_1$ is even and $t_2$ is odd,} \\[2pt]
(\pi,\pi,0), & \textrm{if $t_1$ and $t_2$ are even.}
\end{array}\right.
$$
Then $k_o\in K_1^+=K_2^-$, where $K_1^+=K_2^-$ is defined by \er{sbG3} and for $k_o$ the condition \er{soeq} is fulfilled. For the point $0\in K_1^-=K_2^+$ the condition \er{soeq} is also fulfilled. Thus, due to Corollary \ref{CJig}, the perturbed hexagonal lattice $\cG_t$ is isospectral to its limit graph $\wt\cG$, i.e., $\s(H_{\cG_t})=\s(H_{\wt\cG})$. This and Lemma \ref{Ldeb} yield \er{isE3}. \qq $\Box$

\medskip

\textbf{Acknowledgments.}  I am very grateful to my dear teacher, Professor Evgeny Korotyaev, for his continuous encouragement and inspiration. 



\end{document}